\magnification=\magstep1
\input amstex
\documentstyle{amsppt}

\hsize=5in
\vsize=7.5in

%%%%%%%%%%%%%%%%%%%%%%%%%%%%%%%%%%%%%%%%%%%%%%%%%%%%%%%%%%
\def\pf{\hfill $\blacksquare$}
\def\ub{\underbar}
\def\c{\cite}
\def\fr{\frac}
\def\fg{\frak{g}}
\def\Z{\Bbb Z}
\def\W{\Bbb W}
\def\C{\Bbb C}
\def\call{\Cal L}
\def\ll{\Bbb L}
\def\X{X}
\def\n{\nabla}
\def\uend{\underline{\text{End}}\ }
\def\alt{\text{Alt}}
\def\ext{\text{Ext}}
\def\hom{\text{Hom}}
\def\iff{\text{iff}}
\def\lie{\text{Lie}}
\def\tor{\text{Tor}}
\def\per{\text{per}}
\def\endo{endomorphisms }
\def\hoc{Hochschild }
\def\st{such that }
\def\asa{(A,A)}
\def\wtd{\widehat{D}}

\def\oa{\overline{A}}

\def\ux{\underline{x}}
\def\unn{\underline{\n}}
\def\sd{D\hskip-.09in{/}\,}%\not{D}}
\def\ssd{\not{D}}
\def\usd{\underline{\ssd}}
\def\rc{\overline{C}^{\lambda}}

\def\hc{\overline{HC}}
\def\eH{\Cal{H}}
\def\l{\frak{h}}
\def\glk{\frak{gl}(k)}
\def\bimod{_{\frak{gl}(k[\eta])}U(\frak{gl}(A[\epsilon,\eta]))_{\frak{gl}(k[\epsilon])}}
\def\rightisoarrow{\;\widetilde{\rightarrow}\;}
\def\H{\text{H}}
\def\Ker{\text{Ker}}
\def\Im{\text{Im}}
\def\deg{\text{deg}\;}
\def\ccp{CC_*^{per}}
\def\ah{{\Bbb A}^{\hbar}}
\def\ahc{{\Bbb A}^{\hbar}_c}
\def\E{{\Cal E}_A^*}
\def\hb{\hbar}
\def\*{\star}

\def\Tr{\text{Tr}}
\def\tr{\text{tr}}
\def\Cc{\check{C}}
\def\varksi{\xi}
\def\R{\Bbb R}

\def\Nabla{D}
%%%%%%%%%%%%%%%%%%%%%%%%%%%%%%%%%%%%%%%%%%%%%%%%%%%%%%%%%%
\rightheadtext{  \ \   }
\leftheadtext{  \ \   }
\TagsOnRight
\NoBlackBoxes
\topmatter
\title
On the cohomology ring of an algebra
\\
\endtitle
\author
Ryszard Nest \\
University of Copenhagen \\
and \\
Boris Tsygan\footnote"*"{Partially supported by NSF Grant DMS--9307927 .} \\
Pennsylvania State University
\endauthor
\date
March 1, 1996; revised October, 1997
\enddate
\endtopmatter
\NoBlackBoxes

\TagsOnRight
\document
\baselineskip =1\baselineskip
%\hskip 60mm

To Yu. I. Manin on his 60th birthday
\head
\underbar{\bf Abstract}
\endhead

We define several versions of the cohomology ring of an associative algebra. These ring structures unify some well known operations from homological algebra and differential geometry. They have some formal resemblance with the quantum multiplication on Floer cohomology of free loop spaces. We discuss some examples, as well as applications to index theorems, characteristic classes and deformations.

\subhead
0. \ \ub{Introduction}
\endsubhead

This paper is a continuation and extension of \c{NT2}.

Let $X$ be a smooth manifold or an affine algebraic variety. Let $A$ be the algebra of smooth or regular functions on $X$. It is well known that most geometric properties of $X$ can be recovered from the algebra $A$. By (affine) noncommutative geometry one usually means extending corresponding constructions to the case when $A$ is not necessarily commutative.  

Among differential-geometric objects, the easiest ones to generalize
to the noncommutative setting are vector fields. A noncommutative vector field is just a derivation of an algebra $A.$ As for noncommutative multivector fields and noncommutative differential forms, the answer has been suggested in \c{HRK}. Consider the standard cochain complex $C^*(A,A)$ computing the Hochschild cohomology 
$$
        H^*(A,A)=Ext^*_{A\otimes A^o}(A,A)
$$
(Section 1) and the standard chain complex computing the Hochschild homology
$$
          HH_*(A)=H_*(A,A)=Tor_*^{A\otimes A^o}(A,A)
$$
(Section 2). It was shown in \c{HRK} that, if $A$ is the algebra of regular functions on an affine algebraic variety $X$ over a field $k$ of characteristic zero, then 
$$
     HH_i(A) \simeq \Omega^i_{X/k} 
\tag0.1
$$
$$
      H^i(A,A)\simeq \Gamma (X, \wedge ^{i}TX)
\tag0.2
$$
The standard cochain complex $C^*(A,A)$ is a differential graded associative algebra (with the cup product which can also be interpreted as the Yoneda product). It is also a differential graded Lie algebra (with the Gerstenhaber bracket \c{G}). Under the isomorphism (0.1) these operations turn into the wedge product and the Schouten bracket on multivector fields. The cap product from \c{CE}
$$
   C^p(A,A)\otimes C_n(A,A) \rightarrow C_{n-p}(A,A)
\tag0.3
$$
becomes a noncommutative analogue of contracting a form with a multivector field.

On the other hand, the standard chain complex $C_*(A,A)$ is not a differential graded algebra unless $A$ is commutative. Instead, there is the exterior product \c{CE}
$$
       C_*(A,A)\otimes C_*(B,B) \rightarrow C_*(A\otimes B,A\otimes B)
\tag0.4
$$
which induces an isomorphism on homology (for algebras over a field).

In \c{R} Rinehart defined the differential
$$ B:
C_*(A,A) \rightarrow C_{*+1}(A,A)
$$
which commutes with the \hoc differential and becomes the de Rham
differential under the isomorphism (0.1). A systematic study of a 
noncomutative de Rham theory, or cyclic homology, was undertaken in \c{C},
\c{FT}, \c{L}, \c{T},  etc. One defines the periodic cyclic complex
$(CC_*^{per}(A), B+b) $ of an algebra A, and proves that when $A$ is
the ring of regular functions on a regular affine variety $X$ over a field
$k$ of characteristic zero, the homology of this complex is the de Rham
cohomology of $X$. It is worth mentioning that the correct space of
``noncommutative differential forms'' is not the \hoc homology
$\H_*(A,A)$ but the space of all \hoc chains $C_*(A,A)$.
The \hoc differential  $b$ itself is an odd noncommutative bivector
field and the Hochschild homology is ``a space of physical states''
$\Ker b/\Im b$. We have already mentioned
that there seems to be no way of recovering noncommutative differential
forms together with their multiplicative structure. It turns out, however,
that if one considers all differential operators on $\Omega^*(X)$
(not just  zero order multiplication operators) one gets a ring
which has a noncommutative generalization. The construction we are going
to present rests upon an idea of Yu.I.Manin (\c{M}).

Let $R$ be a monoidal category.  For an object $A$ of
$R$ one calls an object $\uend \ A$ of $R$ the inner object of  \endo
of $A$ if there are natural morphisms
$$
        \uend \ A\times\uend \ A\to\uend \ A\quad;\quad
        \uend \ A\times A\to A
\tag0.5 
$$
which are universal and associative in a natural way.
In the category of associative algebras inner objects of \endo
do not exist (for example, the set of \endo of an algebra is not a 
linear space).  In this paper, however, we will show that, if one takes 
for $\uend \ A$ the differential graded algebra $C^*(A,A)$ of Hochschild
cochains, the maps (0.5), in a sense, still exist.  More precisely,
they exist if one passes from the category of algebras
to the category of complexes by means of some well known homological
functors.

These functors put in correspondence to an algebra its \hoc or periodic 
cyclic complex.  In the case of bar complex the operations analogous
to the ones that we will describe  were discovered by Getzler - Jones and by Gerstenhaber - Voronov (\c{GJ}, \c{GV}).
Therefore one can say informally that in the homotopy category of
(differential graded) algebras the inner endomorphism objects always exist.

The operations arising from the second map (0.5) (in the above sense) play
a crucial role in our approach to index theorems (Section 6); they also allow 
to construct the bivariant JLO cochain and, hopefully, 
a Chern character with values in bivariant entire cyclic cohomology. 
Those operations bear a strong formal resemblance to the
multiplication in the quantum cohomology of symplectic manifolds. In
particular, one can construct ``the Fukaya category'' of an
associative ring; objects of this category are automorphisms of the
ring, and the object $\uend A$ from (0.5) is the object of endomorphisms of the identity. 

The second author recalls Yu.I. Manin's first reaction upon receiving A.Connes' paper on noncommutative geometry in 1982. Commenting on Connes' formula
$$
da = i[F,a]
$$
where $F$ is an odd element, Prof. Manin suggested that perhaps the differentials should belong to the structure sheaf of a more general superscheme, not just of the odd affine line. This idea, together with I.M.Gelfand's notion of an $(\frak{a}, \Cal{D})$-system, is central to this work.

Let us describe the contents of the paper in more detail.

Let $A$ be an associative unital algebra over a commutative unital
ground ring $k$.  Consider $A$ as a bimodule over itself; by $\E$ we denote the differential graded algebra
$(C^*\asa,\delta,\smile)$ which is the standard complex for computing
$\ext_{A\otimes A^{\circ}}^*\asa=H^*\asa$ (the \hoc cohomology)
equipped with the Yoneda product. We recall some well-known properties
of this algebra in Section 1. In Section 2 we construct the maps  of
complexes 
$$
        \bullet : C_*(A,A) \otimes C_*(\E,\E) \to C_*(A,A) 
$$
$$
        \bullet : C_*(\E,\E) \otimes C_*(\E,\E) \to C_*(\E, \E)
\tag0.6
$$
Here $C_*$ stands for the Hochschild complex computing
$HH_*(A)=\tor_*^{A\otimes A^{\circ}}\asa$. In Sect. 3 we construct  similar morphisms
$$
        \bullet : \ccp(A) \otimes \ccp(\E) \to \ccp(A)  
$$
$$
        \bullet : \ccp(\E) \otimes \ccp(\E) \to \ccp(\E)
\tag0.7
$$

where $CC_*^{\per}$ is the periodic cyclic complex. 

In particular, the center $Z_A$ of $A$ is a (differential) subalgebra
of $C^*\asa$.  One gets the usual products
$$
\gathered
        C_*(Z_A)\otimes C_*(A)\to C_*(A)  \\
        CC_*^{\per}(Z_A)\otimes CC_*^{\per}(A)\to CC_*^{\per}(A)
\endgathered
$$
(\c{CE}, \c{HJ}).

On the other hand, let $D$ be a \hoc cochain of $A$ (in particular a
derivation).  Put
$$
        i_D(a)=(-1)^{{\deg a }\cdot { \deg D}} a \bullet D
$$
in the \hoc complex where $D$, being an element of $\E$, is regarded
as a zero-chain of $\E$; let
$$
        I_D(a)=(-1)^{{\deg a} {\cdot} {\deg D}} a \bullet D
$$
but in the periodic cyclic complex; let
$$
        L_D(a)=(-1)^{\deg a(\deg D-1)} a \bullet (1 \otimes D) 
$$

(in this case the formulas are the same for \hoc or periodic cyclic case). 
Then $L_D$ is the Lie derivative and $I_D$ is the contraction operator
of Rinehart \c{R}; the fact that $\bullet$ is a morphism of complexes
implies the Cartan homotopy formula
$$
        [B+b,\;\;I_D]=L_D+I_{\delta D}
\tag0.8
$$
(cf. \c{R}).

There is another product on $C^*\asa$, the Gerstenhaber bracket.  It
makes $C^{*}(A,A)[1]$ a differential graded Lie algebra which we denote by
$(\fg(A),\delta,[\;,\;])$ (cf \c{Ge}).  The Lie derivative $L_D$ turns
$C_*(A)$ and $CC_*^{\per}(A)$ into differential graded
$\fg(A)$-modules.  Let $R_*(\fg(A))=\wedge^*(\fg(A))\otimes U(\fg(A))$
be the standard Koszul resolution of the trivial right $\fg(A)$-module
$k$.  We construct a morphism of complexes of modules (Sect. 4):
$$
        J:R_*(\fg(A))\to\uend CC_*^{\per}(A)
\tag0.9
$$

In our view, the existence of this morphism is one of the main
features of noncommutative differential geometry.  Compare this with
the usual situation: \ let $X$ be a smooth manifold, and let $D_i$ be
vector fields on $X$.  Put
$$
        J(D_1\wedge\dots\wedge D_m)=i_{D_1}\dots i_{D_m}:\Omega_X^*
        \to\Omega_X^{*-m}
$$
Then
$$
\multline
        [d,J(D_1\wedge\dots\wedge D_m)]=\underset{i}\to\sum
        (-1)^{m-i+1}J(\dots\wedge\wtd_i\dots)L_{D_i}+   \\
        +\underset{i<j}\to\sum (-1)^{i+j-1}J([D_i D_j]\wedge\dots\wedge
        \wtd_i\wedge\dots\wedge\wtd_j\dots)
\endmultline
\tag0.10
$$
The morphism (0.9) satisfies a similar formula in the general context
(when, for example, $D_i$ are derivations of $A$; of course, one has 
to replace $d$ by $B+b$).

Note that one could simply put $J(D_1\wedge\dots\wedge
D_m)=I_{D_1}\dots I_{D_m}$ if the operators $I_D$ (anti)-commuted; but
they do not, so the formula is different.

In section 5 we construct another  multiplication which is crucial for our
approach to index theorems. Let $A$ be an associative algebra. By
$\overline{C}_{*-1}^\lambda(A)$ we denote the reduced cyclic complex of $A$;
one can define this complex also as
$$
\overline{C}_{*-1}^\lambda(A)=\text{Prim}\;C_*(\frak{gl}(A),\frak{gl}(k);k)
$$
(the subcomplex of primitive Lie algebra chains; cf. \c{L}). We construct
the operation
$$
\overline{C}_{*-1}^\lambda(A)\otimes CC_*^{per}(A)\longrightarrow CC_*^{per}(A)
$$
or in other words the morphism of complexes
$$
\chi:\ \overline{C}_{*-1}^\lambda(A)\rightarrow\uend CC_*^{per}(A)
$$
To obtain this operation one applies the homomorphism $J$ in the
case when $A$ is replaced by its matrix algebra. Let $M_{\infty}(A)$ be the
algebra of matrices over $A$ with finitely many non-zero diagonals;
let $M(A)$ be its ideal of matrices with finitely many non-zero
entries and let ${\frak{gl}}$ be the same algebra viewed as a Lie algebra. 
We will see in Section 5 that one can modify the operation $J$ to get
a homomorphism 
$$ C^*({\frak{gl}}(A),{ \frak{gl}}(k)) \otimes CC_*^{per}(M_{\infty}(A))\rightarrow CC_*^{per}(M(A))$$

Here ${\frak{gl}}(A)$ is viewed as a subalgebra of the
algebra of inner derivations. To get the operation $\chi$ one uses the
embedding $CC_*^{per}(A) \rightarrow CC_*^{per}(M_{\infty}(A))$ and the trace
map $CC_*^{per}(M(A)) \rightarrow CC_*^{per}(A)$.

The next four sections of the paper are devoted to examples and applications. We study the
 product $\bullet$ on the \hoc complex. One can define the \hoc chain complex
of the algebra of cochains  in two different ways:
$$
\gathered
\Cal{C}^n(A)=\bigoplus_{i+j=n} C_{-i}(\Cal{E}_A^*,\Cal{E}_A^*)^j \\
\Cal{C}^n_\infty(A)=\prod_{i+j=n} C_{-i}(\Cal{E}_A^*,\Cal{E}_A^*)^j
\endgathered
$$
Therefore we get two versions of the cohomology ring, $\eH^*(A)$ and
$\eH^*_\infty(A),$ together with a map $\eH^*\rightarrow \eH^*_\infty$.
First let  $A=k[X]$, where $X$ is regular affine over a ring of characteristic
zero. We show that 
$$
\eH^*(A)=\Cal{D}(\Omega_X^*)^{o};\qquad \eH^*_\infty(A)=\text{End}(\Omega^*_X)
$$
Here $\Cal{D}(\Omega_X^*)$ stands for the ring of differential
operators on the graded module of differential forms, that is to say,
differential operators on the exterior algebra of the cotangent
bundle; the symbol $^{o}$ means the opposite ring).

Next we consider the case of a deformed algebra of functions. Let $(M,\omega)$
be a symplectic manifold and $*$ be a star product on $C^\infty(M)$
(cf. \c{BFFLS}) such that $f*g-g*f=i\hbar\{f,g\}+O(\hbar^2)$.
Isomorphism classes of such star products are parametrized by
$\H^2(M,\Bbb{C}[[\hbar]])$ (\c{LDW}, \c{De}, \c{NT}).
Let $A=\Bbb{A}^\hbar(M)[\hbar ^{-1}]=C^\infty(M)[[\hbar, \hbar^{-1}]$ with the product $*$. Then
$$
\eH^*(A)\rightisoarrow \H^*(M^{S^1},\Bbb{C}[[\hbar, \hbar^{-1}])
$$
if $M$ is simply connected;
$$
\eH^*_\infty(A)\rightisoarrow \H^*(M,\Bbb{C}[[\hbar, \hbar^{-1}])
$$
Note some resemblance with the multiplication on Floer cohomology.

Let us now mention some applications of our constructions.

In Section 6 we outline the proof of the algebraic index theorem from \c{Fe} and \c{NT}. Our approach owes very much to the ideas of B.Feigin and I.M.Gelfand.
Let us mention some other applications, namely to deformations and characteristic classes.

Let $\n$ be an odd element of $\fg(A)$ \st
$$
        \delta\n+\fr12[\n\n]=0.
\tag0.11
$$
($\delta$ is the \hoc differential).  We show that, \ub{formally}, $i.e.$ without regard for convergence, if
$$
        e^{\n}=\underset{n\geq 0}\to\sum
        \fr{\wedge^n\n}{n!}\otimes 1
\tag0.12
$$
in $\wedge^* \fg \otimes U\fg$, then for $X(\n)=J(e^{\n}\otimes 1)$
$$
        [B+b,X(\n)]=L_{\n}X(\n)
\tag0.13
$$

Assume, for example, that $A$ is a $\Z_2$-graded algebra and ${\sd}$
is an odd element of $A$.  Then
$$
        \n=ad(\sd)-\sd^2
$$
is an odd cochain satisfying (0.11).  Using (0.9), we construct a
morphism
$$
        ch(\sd):CC_*^{\per,(0)}(A)\to CC_*^{\per}(A)
$$
where $CC_*^{\per,(0)}$ is the subcomplex of finite cochains:
$$
 CC_n^{\per, (0)}(A)=\underset{i\equiv n(2)}\to\bigoplus
        A\otimes\oa^{\otimes n},
$$
while, as usually,
$$
 CC_n^{\per}(A)=\underset{i\equiv n(2)}\to\prod
        A\otimes\oa^{\otimes n}.
$$
The component
$$
        ch_0(\sd):CC_*^{\per,(0)} (A)\to A
$$
is given by the formula from [JLO], [GS]:
$$
        ch_0(\sd)=\underset{n\geq 0}\to\sum\;\;
        \underset{\Sb t_0+\dots+t_n=1 \\ t_i\geq 0 \endSb}\to\int
        a_0 e^{-t_0\ssd^2}[\sd,a_0]\dots e^{-t_{n-1}\ssd^2}
        [\sd,a_n]e^{-t_n \ssd^2}dt_0\dots dt_{n-1}
$$

Another application: \ let $A$ be an algebra; assume that another
multiplication law on $A$, $*$, is given.  Put
$$
        \n(a,b)=a*b-ab;
$$ 
then $\n$ satisfies (0.11).  The operator $X(\n)$ is a morphism of
complexes
$$
        X(\n): CC_*^{\per,(0)}(A,\cdot)\to CC_*^{\per}(A,*)
$$
When $ab -a*b \in I$ where $I$ is an ideal of $A$ such that $A$ is
$I$-adically complete, then $X(\n)$ gives a well-defined isomorphism
of periodic cyclic complexes, an explicit version of the theorem of
Goodwillie. If one composes this isomorphism with a trace on the
algebra $(A, *),$ one recovers the cocycle of \c{CFS}.
One can expect that, if the multiplications $\cdot$ and $*$ are
``close'', one would be able to use $X(\n)$ to compare Connes' entire cyclic
cohomologies.

In Section 10 we try to clarify the analogy between the multiplication $\bullet$ and the quantum multiplication, as well as to advance B.Feigin's idea that Lagrangian intersections might be related to the cohomology of deformed algebras. We construct ``the Fukaya category'' of an associative algebra $A$; the objects of this category are automorphisms of $A$ and, for two endomorphisms $\alpha$ and $\beta$, the complex $Hom(\alpha, \beta)$ is the twisted cochain complex $C_*(\E, _{\alpha} \Cal{E}^* _{\beta})$. We also construct a functor putting in correspondence to $\alpha$ the twisted chain complex $C_*(A, A_{\alpha})$. Conjecturally, these are an $A_{\infty}$ category and an $A_{\infty}$ functor. When $A$ is a deformed ring of functions on a symplectic manifold then the cohomology of the above complexes is related to their fixed point sets and to the cohomology of certain path and loop spaces.

{\ub{Acknowledgements}} The second author has many fond memories of the time when he was Yu. I. Manin's graduate student in Moscow from 1982 to 1985. Some of the ideas and methods of this paper are due to his collaboration with Yu.L. Daletski, B. Feigin, I.M. Gelfand, as well as to discussions with Yu. I. Manin,
A. Beilinson, V. Schechtman, M. Wodzicki and many others back then
(and later); some other methods and ideas came from the first author's
collaboration with G. Elliot and T. Natsume. We would like to thank
P. Bressler, J.L. Brylinski, A.Connes, B.V. Fedosov, M. Flato, E.Getzler, J. Kaminker, C. Kassel, M. Kontsevich, J.L. Loday, A.Radul, N.Reshetikhin, D. Sternheimer,
D.Tamarkin and A. Voronov for fruitful discussions. Note that
homological operations somewhat similar to ours were first considered by Getzler in \c{G}; there, too, a homological complex of
the differential graded algebra $(C^*\asa,\delta,\smile)$ appeared.
Finally, we would like to thank the referee whose many suggestions
greatly improved the presentation.

\vskip .2in
\noindent\ub{\bf Section 1.}\quad\ub{\bf Hochschild cohomological complex.}
\vskip .2in

Let $A$ be a graded algebra with unit over a commutative unital ring
$k$. Let $M$ be a graded bimodule over $A.$  A \hoc $d$-cochain is a linear map $A^{\otimes d}\to M$.  Put,
for $d\geq 0$,
$$
        C^d (A,M) =\hom_k(\oa^{\otimes d},M)
$$
where $\oa=A/k\cdot 1$.  For a cochain $D$ in $C^*(A,A)$ put
$$
\aligned
        \deg D&=(\text{degree of the linear map }D)+d  \\
        |D|&=\deg D-1; \;\;\;|a|=\deg a - 1
\endaligned
$$
Given a tensor $a_1\otimes\cdots\otimes a_N$ in $A^{\otimes N}$, we
will denote it by $(a_1, \ldots, a_n)$. We will write $\eta_j=\overset{j}\to{\underset{i=1}\to\sum}|a_i|$
(as in \c{G}).  Put for cochains $D$ and $E$ in $C^*(A,A),$
$$
\gathered
        (D\smile E)(a_1,\dots,a_{d+e})=(-1)^{\deg E\cdot\eta_d}
        D(a_1,\dots,a_d)  \\
        \times E(a_{d+1},\dots,a_{d+e});
\endgathered
\tag1.1
$$
$$
\gathered
        (D\circ E)(a_1,\dots,a_{d+e-1})=\underset{j\geq 0}\to\sum
        (-1)^{|E|\eta_j}
        D(a_1,\dots,a_j,  \\
        E(a_{j+1},\dots,a_{j+e}),\dots); \\
        [D, \; E]= D\circ E - (-1)^{|D||E|}E\circ D
\endgathered
\tag1.2
$$
These operations define the graded associative algebra
$(C^*\asa,\deg,\smile)$ and the graded Lie algebra
$\fg ^*(A)=(C^{*+1}\asa$, $|\;\cdot\;|$, [\;,\;]) (cf. \c{CE}; \c{Ge}).
Let
$$
        m(a_1,a_2)=(-1)^{\deg a_1}\;a_1 a_2;
$$
this is a 2-cochain of $A$ (not in $C^2$, because it is a bilinear map
from $A$ to $A$, not from $\overline{A}$ to $A$).  Put
$$
\aligned
        &\delta D=[m,D];  \\
        &(\delta D)(a_1,\dots,a_{d+1})=(-1)^{|a_1||D|+|a_1|+1}\times  \\
        &\times a_1 D(a_2,\dots,a_{d+1})+  \\
        &+\overset{d}\to{\underset{j=1}\to\sum}(-1)^{|D|+\eta_j}
                D(a_1,\dots,a_j\;a_{j+1},\dots,a_{d+1})  \\
        &+(-1)^{|D|+\eta_d+1}D(a_1,\dots,a_d)a_{d+1}
\endaligned
\tag1.3
$$
The last formula defines a differential of degree $+1$ on $C^*(A,M)$ for any $M.$
For an element $x$ of $A$, let $\ux$ be the corresponding zero-cochain in $C^*(A,A).$
By definition
$$
        (\delta\ux)(a)=(-1)^{\deg x}[x,a];
$$
a one-cochain $D$ is a cocycle \iff\ it is a derivation.  One has
$$
\gathered
        \delta^2=0;\quad\delta(D\smile E)=\delta D\smile E+(-1)^{\deg D}
                D\smile\delta E  \\
        \delta[D,E]=[\delta D,E]+(-1)^{|D|}\;[D,\delta E]
\endgathered
$$
($\delta^2=0$ follows from $[m,m]=0$).

Thus $C^*\asa$ becomes a complex; the cohomology of this complex 
is $H^*\asa$ or the \hoc cohomology.  The $\smile$ product induces the
Yoneda product on $H^*\asa=\ext_{A\otimes A^0}^*\asa$.  The operation
[\;,\;] is the Gerstenhaber bracket \c{Ge}. 

If $(A, \;\; \partial)$ is a differential graded algebra then one can define the differential $\partial$ acting on $C^*(A,A)$ by 
$$
\partial D \;\; = \; [\partial , D]
$$
\proclaim
{\ub{Definition 1.1}} Define the differential graded algebra $\E $ as the complex $C^*(A,A)$ with the differential $\delta$ (or $\delta + \partial$ if $(A,\partial )$ is a differential graded algebra), the grading $\text{deg}$ and the product $\smile .$
\endproclaim

 In Sections 1--9 the only cochains we will be considering will be
those from $C^*(A,A)$, {\it{i.e.}} when $A=M.$ For \hoc cochains 
$D_i$ define a new \hoc cochain $D_0\{D_1, \ldots , D_m\}$ by the following formula of Gerstenhaber (\c{Ge}) and Getzler (\c{G}):

$$
D_0\{D_1, \ldots , D_m\}(a_1, \ldots, a_n) = 
$$
$$
=\sum (-1)^{ \sum_{p=1}^{m}\eta_{i_p}| D_p|}  D_0(a_1, \ldots ,a_{i_1} , D_1 (a_{i_1 + 1}, \ldots ),\ldots ,D_m (a_{i_m + 1}, \ldots ) , \ldots)
\tag1.4
$$
The following statements are contained in \c{Ge}, \c{G}, \c{GV}.
\proclaim
{\ub{Proposition 1.2}} One has
$$
\gathered
(D\{E_1, \ldots , E_k \})\{F_1, \ldots, F_l \}=\sum (-1)^{\sum _{q \leq i_p}|E_p||F_q|} \times \\
\times D\{F_1, \ldots , E_1 \{F_{i_1 +1}, \ldots , \} , \ldots ,  E_k \{F_{i_k +1}, \ldots , \}, \ldots, \}
\endgathered
$$
\endproclaim
\demo
{\ub{Proof}}  Direct computation   \pf
\enddemo

\proclaim
{\ub{Example 1.3}}
One has 
$$
D \smile E = (-1)^{\deg D} m\{D,E\} ; \;\;\;\; D\circ E = D\{E\}
$$
\endproclaim
\proclaim
{\ub{Corollary 1.4}} 
$$ 1). \;\;\;\; (D\circ E)\circ F - D\circ (E\circ F) = D\{E,F\}+(-1)^{|E||F|}D\{F,E\};$$
$$ (-1)^{|D||F|}[[D,E],F] + (-1)^{|E||D|}[[E,F],D] + (-1)^{|F||E|}[[F,D],E] = 0 \; .$$

2). Let $M$ be the cochain defined above, but for the
differential graded algebra $\E$: $ M(D_1,D_2)=(-1)^{\deg
D_1}\;D_1 \smile D_2$. For any cochain $D$ in $C^*(A,A)$ define a cochain
${\overline{D}}$ in $C^*(\E, \E)$ by ${\overline {D}}(E_1, \ldots , E_k) = D\{E_1,
\ldots, E_k\} $ for all $k.$ Then

i).
$${\overline {(D \smile E)}} ={\overline {D}} \smile {\overline {E}};$$
ii).
$$
{\overline{\delta D}}= [\delta,{ \overline{D}}]+[M,{\overline{D}}].
$$
\endproclaim
{\bf {Proof.}} To prove the first identity in 1), one applies
Proposition 1.2 in the case when $k=l=1$. The second identity of 1)
follows from the first if one antisymmetrizes the first (with
appropriate signs). To prove the identity 2), i) one applies
Proposition 1.2 to compute $(m\{D,E\})\{F_1, \ldots,F_m\}$ to get the
following formula
$$
 \;\;\;\;(D\smile E)\{F_1, \ldots, F_m\} = \sum_{p=0}^{m}
(-1)^{{\deg E}\sum_{i\leq p}|F_i|} D\{F_1, \ldots, F_p\}\smile
E\{F_{p+1}, \ldots, F_m\} \;\; .
$$
which is equivalent to 2), i).

To prove 2), ii) one applies Proposition 2.2 to compute
$(m\{D\})\{F_1, \ldots ,F_m\},$ $(D\{m\})\{F_1, \ldots , F_m\},$ and
$(D\{F_1, \ldots, F_m\})\{m\};$ their sum, with appropriate sings,
gives the following identity
$$
 \;\;\;\;\delta D(\{F_1, \ldots, F_m\}) = (\delta D)\{F_1, \ldots, F_m\}  -
$$
$$-\sum_{p=1}^{m}(-1)^{\text{deg}D + \sum_{i<p}|F_i|}D\{F_1, \ldots,\delta F_p, \ldots, F_m \} +
$$
$$
+ \sum_{p=1}^{m-1}(-1)^{\text{deg}D + \sum_{i\le p}|F_i|}D\{F_1,
\ldots, F_p \smile F_{p+1} \ldots, F_m \}- 
$$
$$
-(-1)^{\text{deg}D + \sum_{i<m}|F_i|}D\{F_1,
\ldots, F_{m-1} \}\smile F_m +
$$
$$
+ (-1)^ {\text{deg}D |F_1|} F_1 \smile D\{F_2,
\ldots, F_{m} \}
$$

which is equivalent to 2), ii).

\vskip .2in
\noindent\ub{\bf Section 2.}\quad\ub{\bf Operations on the Hochschild Complex.}
\vskip .2in
Let $M$ be a differential graded bimodule over a differential graded
algebra $A$.
Recall that the \hoc homological complex $(C_*(A,M),b)$ is the following:
$$
        C_n(A,M)=M\otimes\oa^{\otimes n}
$$
We will always write, as in Section 1, 
$$
(a_0,\dots,a_n) = a_0 \otimes \dots \otimes a_n
$$
For $a_0$ in $M$ and $a_i$ in $A$, $i > 0,$ define
$$
\multline
        b(a_0,\dots,a_n)=\overset{n-1}\to{\underset{j=0}\to\sum}
        (-1)^{\eta_{j+1}+1}(a_0,\dots,a_j\;a_{j+1},\dots,a_n)  \\
        +(-1)^{(|a_n|+1)(\eta_n+1)+1}
        (a_n a_0,a_1,\dots,a_{n-1}).
\endmultline
\tag2.1
$$
We shall denote $C_*(A,A)$ simply by $C_*(A).$ We introduce a grading on $C_*(A,A)$ by the formula
$$\deg(a_0 , \ldots , a_n) = \sum  \text{deg} a_i + n$$
(A rule for remembering the signs:  let $|\epsilon|=1$, $\epsilon^2=0$;
map $C_n(A)$ to $A[\epsilon]\big/[A[\epsilon],A[\epsilon]]$, $a_0\otimes
\cdots\otimes a_n\mapsto a_0\epsilon a_1\epsilon\cdots a_n\epsilon$;
then $b$ becomes $\frac{\partial}{\partial\epsilon}$).

The homology of $C_*(A, M)$ is the \hoc homology $HH_*(A, M)=\tor_*^{A\otimes
A^{\circ}}(A,M)$. We will denote   $HH_*(A, A)$ simply by $HH_*(A).$ If $A$ is a differential graded algebra and
$\partial$ the differential in $A$, one extends $\partial$ to
$C_*(A)$:
$$
        \partial(a_0,\dots,a_n)=\overset{n}\to{\underset{j=0}\to\sum}
        (-1)^{\eta_j}(a_0,\dots,\partial a_j,\dots,a_n)
\tag2.2
$$
For $a$ in $C_*(A,A)$ and $x$ in $C_*(\E , \E)$ define
$$         
     a\bullet x =  a\bullet _1 x +  a\bullet _2 x
\tag2.3
$$
where
$$
  \multline
     (a_0 ,\ldots , a_n) \bullet_1  (D_0 , \ldots ,      D_m) = \sum
(-1)^{\sum_{i>0}|a_i|\text{deg}D_0 + \sum_{p>0}\sum_{i>i_p}|a_i||D_p|}\times \\
 \times (a_0 D_0(a_1, a_2, \ldots ) ,\ldots ,a_{i_1} , D_1 (a_{i_1 + 1}, \ldots ),\ldots ,D_m (a_{i_m + 1}, \ldots ) , \ldots)
\endmultline
\tag2.4
$$
$$
 \multline
     (a_0 ,\ldots ,a_n) \bullet_2  (D_0 , \ldots ,      D_m) = \\
=\sum_{q \leq n+1} (-1)^{\sum_{i<q}|a_i|\sum_{i\geq q}|a_i| + \sum
_{i_0 < i < q} |a_i| \text{deg}D_0 + \sum _{p>0}\sum
_{i_p < i < q} |a_i| |D_p|} \times \\
\times (D_m(a_q, \ldots , a_n, a_0, \ldots, a_{i_0}) D_0(a_{i_0 + 1}, \ldots ) ,\ldots , a_{i_1} , \\
, D_1 (a_{i_1 + 1}, \ldots ), \ldots , D_{m-1} (a_{i_{m-1} + 1}, \ldots ) , \ldots )
\endmultline
\tag2.5
$$
The sum in (2.5) is taken over all $q, \ i_0, \ldots , \ i_{m-1}$ for which $a_0$ is inside $D_m.$
\proclaim
{\ub{Theorem 2.1}} The map
$$
      \bullet : C_*(A,A) \otimes C_*(\E , \E) \to C_*(A,A)
\tag2.6
$$
is a morphism of complexes, i.e.
$$
       b(a \bullet x) = (ba) \bullet x + (-1)^{\deg a} a \bullet (b + \delta)x
$$
\endproclaim
\demo
{\ub{Proof}} Let $\alpha = (a_0, \ldots, a_n)$ and $D = (D_0, \ldots,
D_m).$ Compute $b(\alpha \bullet _1 D) -( b\alpha) \bullet _1 D. $ It
will consist of the following terms. 

1). All the terms with $a_i a_{i+1}$ or $a_n a_0$ outside $D_j$, $i>0.$
They will cancel out (because they appear once in  $b(\alpha \bullet
_1 D)$ and once in  $(b\alpha)\bullet _1 D$).

1a).$ (a_0 a_1 D_0(\ldots), \ldots D_1(\ldots), \ldots)$

2).  All the terms with $a_i a_{i+1}$ or $a_n a_0$ inside $D_j$,
$i>0.$

3). All the terms with $a_i D_j(\ldots)$ or $D_j(\ldots) a_i,$ $i>0.$

The terms 1a), 2), 3) will cancel out with the term $\alpha \bullet _1
\delta D.$

4).$ (a_0D_0(\ldots)D_1(\ldots), \ldots, D_k(\ldots), \ldots);$
$$(a_0D_0(\ldots), \ldots,D_k(\ldots) D_{k+1}(\ldots), \ldots);$$

These terms will cancel out with all the terms in $\alpha \bullet_1
bD$ except the following:
$$
(a_0D_m(\ldots)D_0(\ldots), \ldots, D_1(\ldots), \ldots)
\tag2.7
$$
5). $(D_m (\ldots)a_0D_0(\ldots), \ldots, D_1(\ldots), \ldots).$

Now compute $b(\alpha \bullet _2 D) -( b\alpha) \bullet _2 D. $
It will consist of the following terms.

6). All the terms with $a_i a_{i+1}$ or $a_n a_0$ outside $D_j.$
They will cancel out (because they appear once in  $b(\alpha \bullet
_2 D)$ and once in  $(b\alpha)\bullet _2 D$).

7).  All the terms with $a_i a_{i+1}$ or $a_n a_0$ inside $D_j.$

8). All the terms with $a_i D_j(\ldots)$ or $D_j(\ldots) a_i.$

The terms 7), 8) will cancel out with some of the terms in  $\alpha \bullet_2 \delta D$

9).$D_m(\ldots a_0 \ldots)D_0(\ldots), \ldots, D_1(\ldots), \ldots).$

9a). $D_m(\ldots a_0 \ldots)D_0(\ldots) D_1(\ldots), \ldots,
D_2(\ldots),\ldots).$

9b).$D_m(\ldots a_0 \ldots)D_0(\ldots), \ldots,
D_k(\ldots)D_{k+1}(\ldots), \ldots).$

The terms 9) and 9a) cancel out with all the terms in $\alpha
\bullet_2 bD$ except the following:
$$
((D_{m-1}\smile D_m)(\ldots, a_0, \ldots)D_0(\ldots), \ldots,
D_1(\ldots), \ldots)
\tag2.8
$$
and 
$$
D_{m-1}(\ldots, a_0, \ldots)D_m(\ldots)D_0(\ldots), \ldots,
D_1(\ldots), \ldots)
\tag2.9
$$
The combination of the terms 9b), (2.8) and (2.9) cancels out. The
remaining terms are (2.7) and 5). But they will cancel out with the
terms in $\alpha \bullet_2 \delta D$ that do not cancel out with 7), 8).

\pf
\enddemo
Consider some examples of the product  $\bullet$. When all $D_i$ are zero-cochains then one gets the shuffle product from \c{CE}. On the other hand, 
put, for a $d$-cochain $D$ and for $a$ in $C_*(A,A)$, $i_D(a) = (-1)^{\deg a\;\deg D}(a \bullet D);$
$$
        i_D(a_0,\ldots,a_n)=(-1)^{\deg a_0 \cdot \deg D_0}
        (a_0 D(a_1,\dots,a_d),\; a_{d+1},\dots,a_n)
\tag2.7
$$

One gets the cap product from \c{CE}:
$$
        C^d(A,A)\otimes C_n(A,A)\to C_{n-d}(A,A)
$$
Now, for a $d$-cochain $D$ and for $a$ in $C_*(A,A)$, let $L_D (a) = (-1)^{|D| \cdot \deg a} a \bullet (1, D);$ one has
$$
\multline
L_D(a_0, \ldots,a_n)= 
\sum_{q=0}^{n-d+2} (-1)^{(\eta _{n+1} - \eta _q) \eta _q + |D|} 
(D(a_q , \ldots, a_0 , \ldots), \ldots, a_{q-1}) +\\ 
+\sum_{k=1}^{n-d}(-1)^{(\eta _{k+1}+1) |D|}(a_0, \ldots, a_k , D(a_{k+1}, \ldots ), \ldots, )
\endmultline
\tag2.10
$$
One has
$$
     [b, L_D] = -L_{\delta D}
$$
Now construct the product

$$
      \bullet : C_*(\E, \E) \otimes C_*(\E , \E) \to C_*(\E , \E)
\tag2.11
$$
as follows. For a cochain $D$ let $D^{(k)}$ be the following $k$-cochain of $\E$:
$$
D^{(k)}(D_1, \ldots, D_k) = D\{D_1, \ldots, D_k\}
$$
\proclaim
{\ub {Proposition 2.2}} The map 
$$
D \mapsto \sum_{k \geq 0} D^{(k)}
$$
is a morphism of differential graded algebras $\E \to \Cal{E}^*_{\E}$.
\endproclaim
\demo
{\ub{Proof}} Follows from  Corollary 1.4, 2).    \pf
\enddemo
Now we construct the product  $\bullet$ by composing the one for the algebra $\E$ with the above map.

The product (2.11) is given by formulas (2.3-2.5) where $a_i$ are now viewed as \hoc cochains and $D(a_1, \ldots, a_d)$ is replaced by $D\{a_1, \ldots, a_d\}$.
\proclaim
{\ub{Proposition 2.3}}
The product (2.11) is homotopically associative.
\endproclaim
We omit the proof.

\vskip .2in
\noindent\ub{\bf Section 3.}\quad\ub{\bf Operations on the periodic cyclic
complex.}
\vskip .2in
In the situation of Sect. 2 define
$$
        B(a_0,\dots,a_n)=\overset n\to{\underset{j=0}\to\sum}
        (-1)^{(\eta_{n+1}-\eta_j)\eta_j}(1,a_j,\dots,a_n, 
        a_0,\dots,a_{j-1})
\tag3.1
$$
Then $B^2=b^2=Bb+bB=0$; denote
$$
        CC_n^{\per}(A)=\underset{i\equiv n(2)}\to\prod
        A\otimes\oa^{\otimes n}
\tag3.2
$$
with the differential $b+B+\partial$ (for the differential graded
algebra $(A,\partial)$).  Put 
$$
\multline
        (a_0, \ldots, a_n) \bullet _3 (D_1, \ldots, D_m) = \\
        =\underset{0\leq p\leq m}\to\sum 
        (-1)^{\underset{r<p}\to\sum (|D_r|+1)\underset{r\geq p}\to\sum                 (|D_r|+1)+\underset{r}\to\sum|D_r|}\times \\
        \times\underset{0\leq j\leq n}\to\sum                           
        (-1)^{(\eta_{n+1}-\eta_j)\eta_j}\lambda'
        (D_p,\dots,D_m,D_0,\dots,D_{p-1})  
        (1,a_j,\dots,a_n,\;a_0,\dots,a_{j-1})
\endmultline
\tag3.3
$$
Here

$$
\multline
        \lambda(D_1,\dots,D_m)(a_0,\dots,a_n)=  \\
        =\underset{j_1>0}\to\sum
        (-1)^{|D_1|(\eta_{n+1}-\eta_{j_1})+\dots+|D_m|(\eta_{n+1}-\eta_{j_m})}\times \\
        \times(a_0,\dots,D_1(a_{j_1}\dots),\dots,D_m(a_{j_m},\dots),\dots);  \\
\endmultline
\tag3.4
$$

 $\lambda'$ is the sum of all those terms in $\lambda$ 
which contain $D_0(\dots)$ after $a_0$.

\proclaim
{\ub{Theorem 3.1}} The map
$$
      \bullet : \ccp(A) \otimes \ccp(\E) \to \ccp(A)
\tag3.5
$$
is a morphism of complexes.
\endproclaim
\demo
{\ub{Proof}} Let $ \alpha = (a_0, \ldots, a_n)$ and $D = (D_0, \ldots,
D_m).$
Compute $b(\alpha\bullet _3 D) - (b \alpha)\bullet _3 D - 
\alpha\bullet _3 bD.$ It will consist of the following terms.

1). All the terms with $a_i a_{i+1}$ or $a_n a_0$ outside $D_j$.
They will cancel out (because they appear once in  $b(\alpha \bullet
_3 D)$ and once in  $(b\alpha)\bullet _3 D$).

2).  All the terms with $a_i a_{i+1}$ or $a_n a_0$ inside $D_j$. 

3). All the terms with $a_i D_j(\ldots)$ or $D_j(\ldots) a_i.$

The terms 2), 3) except for the following one (formula (3.6) will cancel out with the term $\alpha \bullet _3
\delta D.$
$$
(1, \ldots, a_0D_0(\ldots), \ldots, D_1(\ldots), \ldots)
\tag3.6
$$
4). All the terms with $D_i(\ldots)D_{i+1}(\ldots)$ or
$D_m(\ldots)D_0(\ldots)$ but for the following (3.7) will cancel out.
$$
(1, \ldots, D_m(\ldots, a_0, \ldots)D_0(\ldots), \ldots, D_1(\ldots),
\ldots)
\tag3.7
$$
5).$ (D_k(\ldots, a_0, \ldots), \ldots, D_0(\ldots), \ldots).$

6).$ (a_0, \ldots, D_0(\ldots), \ldots).$

7). $D_0(\ldots), \ldots, a_0, \ldots, D_m(\ldots), \ldots).$

(Here $a_0$ may be inside or outside $D_i$).

8). $(a_k, \ldots, a_0, \ldots, D_0(\ldots), \ldots)$

where $a_0$ may be inside or outside $D_i$. These terms will cancel
out because they will enter twice (namely, in the first and in the
last summands of $b(\alpha \bullet _3 D)$).

The term 5) will cancel out with $\alpha \bullet_2 BD,$ the term 6)
with $\alpha \bullet_1 BD,$ and the term 7) with $B\alpha \bullet _1
D.$ The term (3.6) cancels out with $B(\alpha \bullet _1 D)$ and the
term (3.7) with $B(\alpha \bullet_2 D).$ Now, $B\alpha \bullet_2 D=0$
and $B\alpha \bullet_3 D=\alpha \bullet _3 BD = B(\alpha \bullet_3
D)=0,$ which finishes the proof.
  \pf
\enddemo
For a cochain $D$ and for $a$ in $\ccp(A)$, let $I_D(a) = (-1)^{\deg D \cdot \deg a}(a \bullet D);$ $L_D (a) = (-1)^{|D| \cdot \deg a} a \bullet (1, D);$ one has

$$
        [B+b,\;I_D]-I_{\delta D}=L_D.   
$$
\pf

This is the homotopy formula of Rinehart \c{R}.

One can define the product
$$
      \bullet : \ccp(\E) \otimes \ccp(\E) \to \ccp(\E)
\tag3.8
$$
exactly in the same way as in the end of Section 2.

\proclaim
{\ub{Remark 3.2}} Considering $CC^{per}_*(A)=CC^{per}_*(\Cal{E}_A^0)$ as a part of $CC^{per}_*(\E)$ and restricting the product $\bullet$ to it, we get the Hood-Jones product on $CC^{per}_*(A).$ When $A$ is commutative then $CC^{per}_*(A)$ is a subcomplex of $CC^{per}_*(\E)$ and one gets a product on the complex $CC^{per}_*(A).$
\endproclaim
\vskip .2in
\noindent\ub{\bf Section 4.}\quad\ub{\bf The Lie algebra complex.}
\vskip .2in

Let $(\fg,\partial)$ be a differential graded Lie algebra.  Put
$\wedge^*\fg=T(\fg)/<D\otimes E-(-1)^{(|D|-1)(|E|-1)}E\otimes D>$;
$$
\gathered
        R_*(\fg)=\wedge^*\fg\otimes U(\fg);  \\
        \partial^{\text{Lie}}(D_1\wedge\dots\wedge D_m\otimes f)=
        \underset{i}\to\sum(-1)^{\theta_i}D_1\wedge\dots  \\
        \wedge\wtd_i\wedge\dots\wedge D_m\otimes D_i f+  \\
        +\underset{i<j}\to\sum(-1)^{\theta_{ij}}D_1\wedge\dots
        \wedge\wtd_i\wedge\dots\wedge[D_i D_j]\wedge\dots\otimes f
\endgathered
$$
where 
$$
\align
        \theta_i&=\underset{r<i}\to\sum(|D_r|+1)+|D_i|
        \underset{r>i}\to\sum (|D_r|+1);  \\
        \theta_{ij}&=\underset{r<i}\to\sum(|D_r|+1)+|D_i|
        (\underset{i<r<j}\to\sum (|D_r|+1)+1);  
\endalign
$$
Put also
$$
        \partial(D_1\wedge\cdots\wedge D_m)=\sum
        (-1)^{\mu_j}(D_1,\dots,\delta D_j,\dots,D_m)
$$
where
$$
        \mu_j=1+\underset{r<j}\to\sum(|D_r|+1).
$$
\underbar{A rule for remembering the signs}:  Let
$|\epsilon|=1$, $\epsilon^2=0$; map $R_*(\fg)$ to 
$U(\fg[\epsilon])$ by $(D_1\wedge\dots\wedge D_m)\otimes f\mapsto
\epsilon D_1\dots \epsilon D_m\cdot f$; then $\partial^{\lie}$ becomes
$\fr{\partial}{\partial\epsilon}$; $\partial$ is the induced derivation on $U$ but with the opposite sign.

Recall that in Section 3 we defined Lie derivatives $L_D: CC^{per}_*(A) \mapsto CC^{per}_*(A)$ for $D$ in $C^*(A,A).$

\proclaim
{\ub{Lemma 4.1}}
$$[L_D,L_E]=L_{[D,E]}$$
\endproclaim
{\ub{Proof}}. Straightforward.       \pf

Put $\fg(A)=(C^*\asa;\;|\;\cdot\;|;\;[\;,\;];\;\delta)$.  One sees that
$\fg(A)$ is a differential graded Lie algebra and that $U(\fg(A))$
acts on the complex of endomorphisms $\uend CC_*^{\per}(A)$ from the right: $X\mapsto X\cdot L_D$
for $D\in\fg(A)$.  For $f\in U(\fg(A))$, we will denote the
corresponding operator by $L_f$.

The purpose of this Section is to show that the complex $R_*(\fg (A))$ acts $\fg(A)$ - equivariantly on the complex $CC^{per}_*(A)$ , thus generalizing the standard calculus from the commutative case.

Let us consider the complex $CC^{per}_*(\E)$ with the Hood-Jones product (Remark 3.2). In other words, we view $\E$ as $C_0(\E)$ and consider the product $\bullet$ from section 3. We denote this product by $\*.$
For example, 
$$(D) \star (E) = \pm (D \smile E) \pm (1, D, E) \pm (1, E, D)$$
while 
$$(D) \bullet (E) = \pm (D \smile E) \pm (1, D, E) \pm (1, E, D) \pm D\{E\}$$

\proclaim
{\ub{Definition 4.2}}
Put for $\alpha$ in $CC^{per}(A)$ and $D_i$ in ${\fg(A)}$
$$
\gathered
J(D_1\wedge\dots\wedge D_m)\alpha =(-1)^{\deg \alpha\sum_{i=1}^{m}\text{deg}D_i}\fr{1}{m!}\underset{\sigma}\to\sum
\epsilon_{\sigma} \alpha \bullet (D_{\sigma _1} \* (D_{\sigma _2} \* (\ldots \* D_{\sigma _m}))\ldots )
\endgathered
$$
where the sign $\epsilon_{\sigma}$ is taken with respect to the rule
$$
        \epsilon_{(i,i+1)}=(-1)^{(|D_i|-1)(|D_{i+1}|-1)}
$$
for any transposition $(i,\; i+1).$

If $f$ is in $U(\fg(A))$ put
$$
        J(D_1\wedge\dots\wedge D_m\otimes f)=
        J(D_1\wedge\dots\wedge D_m)L_f.
$$

\endproclaim

\proclaim
{\ub{Theorem 4.3}} \roster
\item"{i)\;\;}"  The map $J: R_*(\fg(A)) \to \uend ( CC_*^{per}(A))$ is a homomorphism of right $\fg(A)$-modules;
\item"{ii)\;\;}"  for $\gamma$ in $R(\fg(A)$
$$[B+b,\;J(\gamma)]=J((\partial^{\text{Lie}}+\delta)
        \gamma)$$\item"{iii)\;\;}" If $E_1, \ldots, E_n, D$ are in $\Cal{E}^0$ or in $\Cal{E}^1$ then
$$ [L_D, J(E_1 \wedge \ldots \wedge E_n)] = \sum_{j=1}^n (-1)^{|D|\sum_{i<j}(|E_i|+1)}J(E_1 \wedge \ldots [D,E_j] \wedge \ldots E_n)
$$
\endroster
\endproclaim

\demo
{\ub{Proof}}  i) is obvious.  To prove ii) and iii), one needs the following partial associativity properties of the product $\bullet.$
\proclaim
{\ub{Lemma 4.4}} Let $\alpha \in CC^{per}_*(A)$, $\beta \in   CC^{per}_*(\E)$ and $\gamma \in \E.$ Then
\roster
\item"{i)\;\;}" $\alpha \bullet ((1,D) \bullet \beta) = (\alpha \bullet (1,D)) \bullet \beta;$
\item"{ii)\;\;}" $(\alpha \bullet \beta) \bullet (1,D) = \alpha \bullet (\beta \bullet (1,D))$ if $D$ is in $\Cal{E}^0$ or in $\Cal{E}^1;$ 
\item"{iii)\;\;}"$$\alpha \bullet (1,D) \bullet \beta  - (-1)^{|\beta||D|} \alpha \bullet \beta \bullet (1,D) = \alpha \bullet [D, \beta]$$ where $D$ is in $\Cal{E}^0$ or in $\Cal{E}^1$ and $[D, \beta]$ stands for the usual action of $\fg(A)$ on its tensor powers;
\item"{iv)\;\;}" 
$$
\gathered
(1,D) \bullet (E_0, \ldots, E_n) = (1,D) \* (E_0, \ldots, E_n) + \\
+\sum _{j=0}^n (-1)^{|D|(\sum_{i \leq j}|E_i| + 1)} (E_0, \ldots, E_j \circ D, \ldots, E_n)
\endgathered
$$
\item"{v)\;\;}" $D \* ((1,E) \* \gamma) = (-1)^{(|D|-1)|E|} (1,E) \* (D \* \gamma) $
\endroster
\endproclaim
\demo
{\ub {Proof.}} Direct computation. \pf
\enddemo

Now let us prove the theorem. We shall omit the signs for simplicity. One has
$$
\gathered
[B+b, J(D_1\wedge\dots\wedge D_m)]\alpha = \alpha \bullet \text{Alt}(B+b+\delta)((D_1 \* (D_2 \* (\ldots \* D_m))\ldots )) =\\
= \alpha \bullet (\text{Alt} \sum_j \pm (D_1 \* \ldots ((1, D_j) \* (\ldots \* D_m))\ldots ) + \\
\sum_j \pm (D_1 \* \ldots (\text{ad}( D_j) \* (\ldots \* D_m))\ldots ) + \\
+ \sum_j \pm (D_1 \* \ldots (\delta(D_j) \* (\ldots \* D_m))\ldots ))
\endgathered
$$  
The second sum contributes terms $D_i \smile D_j - (-1)^{(|D_i|-1)(|D_j|-1)}D_j \smile D_i$ and therefore vanishes under antisymmetrization. The third sum contributes the term $J(\delta \gamma)\alpha.$ The first sum, because of (v), is equal to 
$$
\multline
 \alpha \bullet \text{Alt} \sum_j (1,D_j) \*(\pm (D_1 \* \ldots (\widehat{D_j} \* (\ldots \* D_m))\ldots )= \\
 = \alpha \bullet \text{Alt} \sum_j \pm  (1,D_j) \bullet(D_1 \* \ldots (\widehat{D_j} \* (\ldots \* D_m))\ldots )+ \\
+ \alpha \bullet \text{Alt}   \sum_j \pm (D_1 \* \ldots(D_i \circ D_j \* (\ldots (\widehat{D_j} \* (\ldots \* D_m))\ldots )+ \\
= \sum_j \pm J(D_1 \wedge \ldots \wedge \widehat{D_j} \wedge \ldots \wedge D_m)(L_D \alpha) +   \\
+\sum_{i<j} \pm J([D_i, D_j] \wedge \ldots \wedge \widehat{D_i} \wedge \ldots \widehat{D_j}\wedge \ldots \wedge D_m)( \alpha)
\endmultline
$$
\pf
\enddemo

\vskip .2in
\noindent\ub{\bf Section 5.}\quad\ub{\bf The characteristic map $\chi$}

In this section we construct a natural morphism of complexes
$$
\overline{C}_{*-1}^\lambda(A)\otimes CC_*^{per}(A)\rightarrow CC_*^{per}(A)
$$
where  $\overline{C}_{*-1}^\lambda(A)$ is the reduced cyclic complex equipped with the differential $b$:
$\overline{C}_n^\lambda(A)=\overline{A}^{\otimes(n+1)}/\Im(1-\tau)$
where $\overline{A}=A/k\cdot 1$ and
$$\tau(a_0, \dots, a_n)=
(-1)^{\eta_n(\eta_n-\eta_{n-1})}(a_n, a_0,\dots, a_{n-1}).$$
Note that this complex is well defined (\c{L}).

Let $M(A)$ be the algebra of matrices
$(a_{ij})_{1\leq i, j\le\infty}$ for which $ a_{ij}\in A$ and all but finitely many of
$a_{ij}$ are zero. The same space considered as a Lie algebra is
$\frak{gl}(A)$. The formula 
$$
(a_0,\dots\, a_n) \mapsto (-1)^{\sum \deg a_i}
E_{01}^{a_0}\wedge E_{12}^{a_1}\wedge\dots\wedge E_{n0}^{a_n}
$$
defines a morphism of complexes
$$
\overline{C}_{*-1}^\lambda(A)\ \rightarrow\ C_{*}(\frak{gl}(A),\frak{gl}(k);k)
$$
(with the relative Lie algebra chains in the right hand side). This map
is an isomorphism of $\overline{C}_*^\lambda(A)$ with  the subcomplex
of primitive elements; cf. \c{T},\c{FT}, \c{LQ}.

Let ${M}_\infty(A)$ be the associative algebra of matrices
$(a_{ij})_{1\leq i, j\le\infty}$ such that $a_{ij}\in A$ and $a_{ij}=0$
on all but finitely many diagonals. Consider the differential graded
Lie subalgebra $C^{0} (M(A),M(A)) + Z^{1} (M(A),M(A)) $ There is a morphism to this
algebra from 
the differential graded Lie algebra
$({\frak{gl}}(A[\eta])$ where $\eta$ is a formal
parameter of degree $-1$ such that $\eta^2=0$; the differential
on $A[\eta]$ is $\partial/\partial\eta$. This morphism sends $\eta a$
to $a$ in $C^0$ and $a$ to $\text {ad}(a)$ in $C^1$. From Theorem 4.3 one gets the pairing
$$
{R}_*(\frak{gl}(A[\eta]))\otimes CC_*^{per}({M}_\infty(A))\longrightarrow
CC_*^{per}(M(A))
\tag5.1
$$

Consider the following morphisms:
$$
\gather
CC_*^{per}(A) \buildrel i\over\longrightarrow CC_*^{per}({M}_\infty(A));\\
CC_*^{per}({M}(A))\buildrel tr\over\longrightarrow CC_*^{per}(A)  
\endgather
$$

The map $i$ is induced by the inclusion $a\mapsto a\cdot 1$; the second
map $tr$ acts as follows:
$$
(a_0 m_0,\dots, a_n m_n) \longmapsto
tr(m_0\dots m_n)\; (a_0,\dots, a_n)
$$
for $a_i\in A$ and $m_i\in {M}(k)$. Composing (5.1) with
these maps we obtain
$$
{R_*}(\frak{gl}(A[\eta]))\otimes CC_*^{per}({M}_\infty(A))\longrightarrow
CC_*^{per}(A)
$$
or
$$
U(\frak{gl}(A[\epsilon,\eta])) \longrightarrow \text{End}(CC_*^{per}(A))
\tag5.2
$$
where $[\epsilon,\eta]=[\epsilon,\epsilon]=[\eta,\eta]=0$,
$\deg\epsilon=\deg\eta=-1$ and the differential in the L.H.S.
is $\dfrac\partial{\partial\epsilon}+\dfrac\partial{\partial\eta}$. For modules over any Lie algebra $\fg$ we shall write $\otimes _{\fg}$ instead of $\otimes _{U({\fg})}.$
\proclaim
{\ub{Lemma 5.1}} For any Lie Algebra $\frak{g}$ and any right (left) module
$M$ the  map (5.2) descends to the morphism
$$
k \otimes \bimod \otimes k \longrightarrow \uend ( CC_*^{per}(A))
$$
\endproclaim
{\ub{Proof}} First it is easy to see that the map (5.2) descends to
$k \otimes _{\frak{gl}(k[\eta])}U({\fg}[\epsilon, \eta])$. In fact, this subalgebra acts by the operators
$L_m$ and $L_{\text{ad}(m)}$ where $m\in\frak{gl}(k)$, both are
zero when followed by the trace map $\text{tr}$. To see that one can pass to
$k \otimes _{\frak{gl}(k[\eta])}U(\frak{gl}(A[\epsilon,\eta]))_{\frak{gl}(k[\epsilon])}\otimes k$
note that $L_{\text{ad}(m)}$ is zero on the image of $i$ and that, because
of the explicit formula in Definition 4.2, the operator
$J(\epsilon D_1\cdot\dots\cdot\epsilon D_m)$ is zero on the image of $i$ if
at least one of the $D_i$ is of the form $\text{ad}(m), m\in\frak{gl}(k)$.
Indeed, this formula involves terms $D_i(a_j)$ for any $D_i$ which
is a one-cochain. \pf

Finally there is an explicit map
$$
C_*(\frak{gl}(A[\eta]),\frak{gl}(k);k) \longrightarrow k \otimes \bimod \otimes k
\tag5.3
$$
Let us construct this map. Let $\frak{g}=\frak{gl}(A)$ and $\l=\glk$.
To describe the image of a chain
$D_1\wedge\dots\wedge D_n \wedge \eta E_1\wedge\dots\wedge\eta E_m$,
write the expression
$$
D_1(\epsilon-\eta)\cdot\dots\cdot D_n(\epsilon-\eta)\cdot E_1\epsilon\eta
\cdot\dots\cdot E_m\epsilon\eta
$$
and then represent it as a sum
$$
\sum\pm D_{i_1}\eta\cdot\dots\cdot D_{i_k}\eta\cdot
E_1\epsilon\eta\cdot\dots\cdot E_n\epsilon\eta\cdot
D_{j_1}\epsilon\cdot\dots\cdot D_{j_m}\epsilon
$$
in the symmetric algebra of the graded space $\frak{g}[\epsilon,\eta]$.
For example $D\mapsto D\epsilon-D\eta$;
$$
\gather
D_1\wedge D_2\longmapsto
D_1\eta\cdot D_2\eta+D_1\eta\cdot D_2\epsilon+
(-1)^{(|D_1|+1)(|D_2|+1)}D_2\eta\cdot D_1\epsilon+
D_1\epsilon\cdot D_2\epsilon \\
D_1\wedge D_2\eta\longmapsto-D_1\eta\cdot D_2\epsilon\eta+
(-1)^{(|D_1|+1)(|D_2|+1)}D_2\epsilon\eta\cdot D_1\eta
\endgather
$$
etc.

It turns out that this yields the morphism of complexes
$$
C_*(\frak{g}[\eta],\l;k)\longrightarrow k \otimes \bimod \otimes k
\tag5.4
$$
or, which is the same,
$$
k \otimes _{\frak{gl}(A[\eta])} U(\frak{gl}(A[\epsilon, \eta])_{\frak{gl}(k[ \epsilon])} \otimes k 
\longrightarrow k \otimes \bimod \otimes k
$$
\remark
{\ub{Remark }} These maps are quasi-isomorphisms; one can easily show that
the cohomology of both complexes is
$$
C_*(\frak{g}[\eta],\l;k)=S(\l)_\l=H_*(\text{BH},k)
$$
where $H$ is the Lie group with the Lie algebra ${\frak{h}}$ and $BH$ is
the classifying space of $H.$
\endremark

At the level of cohomology the inclusion
$$
C_*(\frak{g},\l;k)\longrightarrow
C_*(\frak{g}[\eta],\l;k)
\tag5.5
$$
is the dual map to the Chern-Weil map
$$
C^*(\frak{g},\l;k)\longleftarrow W^*(\l)_{basic};
$$
$W^*(\l)=C_*(\l[\eta],k)$ is the Weil algebra of $\l$.
Another realization of the map which is induced by 5.5) on homology is
$$
k \otimes_{\l[\epsilon]}U(\frak{g}[\eta])\otimes_{U(\frak{g})}U(\frak{g}[\epsilon])_{\l[\epsilon]}\otimes k
\longrightarrow
k \otimes _{\l[\eta]}U(\frak{g}[\epsilon,\eta])_{\l[\epsilon]}  \otimes k
\tag5.6
$$

\vskip .2in
\noindent\ub{\bf Section 6.}\quad\ub{\bf Index theorems}
\vskip .2in

To illustrate our methods, we will sketch the proof of the algebraic index theorem from \c{Fe}, \c{NT}. This proof, unlike the original ones, is designed to work in a much more general situation (families, foliations, manifolds with boundaries or with corners, complex varieties, $D$-modules, etc.).

Let $(M, \omega)$ be a $2n$-dimensional symplectic manifold with a star product $* .$
$$
   f *  g = \sum (i \hb)^k \varphi _k (f,g).
$$
Let $\ah (M) = C^{\infty}(M)[[\hbar]]$  with the product $*$. Put also  $\ah _c (M) = C_c^{\infty}(M)[[\hbar]].$ The following is proved in \c{Fe}, \c{NT}. Consider the canonical trace
$$
     \Tr : \ahc (M) \rightarrow \C[\hb ^{-1} , \hb]]
$$
constructed in \c{Fe}, \c{NT}. One has
$$
     \Tr (f) = \frac{1}{(i \hb) ^n n!} [\int_M {f \omega ^n + \sum _{k=1}^{\infty}{(i \hb)^k D_k (f) \omega ^n }}]
$$
where $D_k$ are certain differential operators. All traces on $\ahc (M)$ are proportional.

Let $P^2 = P$, $Q^2 = Q$ be two elements of the matrix algebra $M_N(\ah (M))$ such that $P - Q$ is in $M_N(\ahc (M))$. Let $P_0 = P \text{mod} \hb$,  $Q_0 = Q \text{mod} \hb$, $\text{ch} (P_0) = \sum_k {\frac{1}{k!}\tr P_0 (dP_0)^{2k}}$
(the Chern character of the connection $d P_0 d$ in the vector bundle
$P_0 \C ^N$). We regard $TM$ as a $U(n)$-bundle (since $U(n)$ is a
maximal compact subgroup of $Sp(2n)$). 

Recall that for each deformation quantization of $M$ a characteristic
class $\theta$ is defined which is an element of
$H^2(M)[[\hbar]]$: the class $\theta /i\hbar$ is the curvature of
a Fedosov connection defining the deformation (cf. \c{Fe}, \c{NT}).

\proclaim
{\ub{Theorem 6.1}} 
$$
     \Tr(P - Q) = \int _M { (\text{ch}(P_0) -  \text{ch}(Q_0)) \widehat{A}(TM)e^{\theta / i \hb}}
$$
where $\theta$ is the characteristic class of $*$.
\endproclaim

\demo
{\ub{Sketch of the proof}}  Consider the characteristic map
$$
\chi : \rc _{*-1} (\ah (M))[\hb ^{-1}] \rightarrow \uend \ccp (\ahc (M))[\hb ^{-1}]
$$
Note that the right hand side is the space of global sections of a
fine sheaf ( the left hand side is not). Therefore one can extend
$\chi$ to the morphism to the Chech - cyclic double complex:
$$
\chi : \Cc ^{-*}(\rc _{* - 1} (\ah (M)))[\hb ^{-1}] \rightarrow \uend \ccp (\ahc (M))[\hb ^{-1}]
\tag6.1
$$
Next step is to construct the fundamental class $\Omega$ in the left hand side of (6.1). Note that for any contractible Darboux chart $U_0$ 
$$
\hc_{2n-1} (\ah (U_0))[\hb ^ {-1}] \rightisoarrow \C [ \hb^{-1}, \hb ]]
\tag6.2
$$
(the canonical generator is represented by the cycle 
$$ \Omega _0 = \frac{1}{2n (i \hb)^n} \alt (\varksi _1 \otimes x_1 \ldots \varksi _n \otimes x_n)
\tag6.2.1
$$
 where alternation is over the group $\Sigma _{2n}$ and $[\varksi _i , x_j] = i \hb \delta _i^j$);
$$
\hc_{i} (\ah (U_0))[\hb ^ {-1}] = 0
\tag6.3
$$
for $i > 2n-1.$

Using the spectral sequence $E^2_{ij}=H^{-i}(M, \hc _j (\ah ))$ converging to $\hc _{j-i} (\ah (M)),$ we
construct a cycle of the complex $\Cc ^{-*}(\rc _{* - 1} (\ah (M))[\hb
^{-1}])$ whose homology class, being restricted to any contractible
Darboux chart, gives the canonical generator of $\hc_{2n-1}(\ah (M)).$ Let
$\alpha$ be a periodic cyclic chain of $ \ahc (M); $ $\alpha _0 =
\alpha \text{mod} \hb$ is a periodic cyclic chain of $C^{\infty}_c
(M).$ (For non-unital algebras the periodic cyclic complex is defined
as the kernel $\text {Ker} (CC^{per}(\tilde{A}) \rightarrow CC
^{per}(k))$ were $\tilde{A}$ is the algebra $A$ with adjoined unit). One proves that 
$$
     \Tr (\chi (\Omega))(\alpha) = \int _M \mu (\alpha _0) \text{mod} \hb 
\tag 6.4
$$
where $\mu : \ccp (C^{\infty}_c (M)) \rightarrow (\Omega^*(M), d)$ is the morphism of complexes given by 
$$
   \mu (a_0 \otimes \ldots \otimes a_k) = \frac {(-1)^k}{k!} a_0da_1 \ldots da_k
$$
(\c{L}). Formula (6.4) follows easily from (6.2.1) and from the fact that $[f,g]=i \hb \{ f, g \} \text{mod} \hb ^2 .$

Now consider the cycles $(k-1)! \eta ^{\otimes k}$ in $\Cc ^{-*}(\rc _{* - 1} (\ah (M)[\eta]))[\hb ^{-1}].$ One has
$$
\chi ((k-1)! \eta ^{\otimes k}) = S^k
\tag6.5
$$
where $S$ is the Bott isomorphism.
We know from (6.4) the composition of the canonical trace with
$\chi(\Omega).$ To compute the canonical trace itself we have to
compare the map $\chi(\Omega)$ to the identity map. In view of (6.5),
it is enough to express $\Omega$ in terms of the classes $(k-1)! \eta
^{\otimes k}.$
\proclaim
{\ub {Theorem 6.2}}
In $\Cc ^{-*}(\rc _{* - 1} (\ah (M)[\eta])[\hb ^{-1}]$ $\Omega$ is cohomologous to
$$
 \sum_{k} (-1)^k (\widehat{A}(TM) e ^{\theta / i \hb })^{-1}_{2k} \cdot (n+k-1)! \eta ^{\otimes (n+k)}
$$
\endproclaim

Note that, since $\uend (\ccp (\ahc (M)))$ is the space of global sections of a fine sheaf, its cohomology is a module over the algebra $H^* (M, \C):$ if $(c_{U_0 U_1 \ldots U_p})$ is in $\Cc ^p (M, \C)$   then $\uend (\ccp (\ahc (M)))$  
$$
    c \cdot \alpha = \sum c_{U_0 U_1 \ldots U_p} \rho _{U_0} [B+b, \rho _{U_1}] \ldots [B+b, \rho _{U_p}]
$$
where $\{ \rho _U \}$ is a partition of unity. 

We see that for any periodic cyclic cycle $\alpha$ of $\ccp (\ahc (M))$

$$
    \Tr (e^{-\theta / i \hb}\cdot \alpha) = \int _M {\mu (\alpha _0) \widehat{A}(TM) \text{mod}( \hb)}
\tag6.6
$$
Indeed,
$$
      \Tr (e^{-\theta / i \hb}\cdot \alpha) = \Tr\sum_{k} (-1)^k (\widehat{A}(TM) e ^{\theta / i \hb })^{2k} \cdot S^{n+k} \chi (\Omega)(e^{-\theta / i \hb}\cdot \alpha)
$$
by Theorem 6.2; now we use (6.4).

We apply (6.6) to $\alpha = \text{ch}(P) - \text{ch}(Q)$ where $\text{ch}$ is the Connes-Karoubi Chern character cycle in the periodic cyclic complex:
$$
    \text{ch}(P) = \tr (P + \sum_{k \geq 1} (-1)^k  \frac{(2k)!}{k!}(P - \frac {1}{2}) \otimes P^{\otimes 2k})
\tag6.7
$$
(\c{L}). Note that $\mu (\text{ch}(P)) = \text{ch}(P_0).$
The left hand side of the last formula stands for the Chern character
of the vector bundle $\text {Im} P_0.$

It remains to show that for a periodic cyclic cycle $\alpha$ of $\ah (M)$ {\it over $\C$} (i.e when all the tensor products are taken over $\C$, not $\C [[\hb]] $) one has
$$
\frac{\partial}{\partial \hb} \Tr (e^{-\theta / i \hb}\cdot \alpha) = 0
\tag6.8
$$
This can be done by passing to a one-parameter family of deformations
$$
   f * _{\lambda} g = \sum (i \hb \lambda)^k \varphi _k (f,g).
$$
One gets the deformed algebra $\ah (M \times \R ^1)$. 

It turns out that for any deformed family $M \rightarrow B$ of symplectic manifolds one can construct the canonical morphism of complexes (\c{NT})
$$
   \Tr _{\theta} : \ccp (\ahc (M)) \rightarrow \Omega _c ^* (B)[\hb ^{-1}, \hb]]
$$
such that the $\Omega ^0$ component is, for $b \in B$, 
$$
\Tr _{\theta} (\alpha)(b) = \Tr (\alpha (b) \cdot e^{- \theta _b / i \hb})
$$
where $\alpha (b)$ is the restriction of $\alpha$ to $\ahc (M_b)$ and $\theta _b $ is the characteristic class of $\ah (M_b)$.

It remains to say a few words about the proof of Theorem 6.2. To prove it one has to construct and to compare various cocycles (cochains) of the complex $\Cc ^* (M, \call ^*)$ where $\call ^*$ is the sheaf of complexes $\rc _{-*}(\ah (M))[\hb ^{-1}]$. If $\call ^*$ were a constant sheaf then one would be able to construct such cochains, for example the ones representing characteristic classes of the tangent bundle. One observes that $\call ^*$ is {\it constant up to homotopy} in the following sense. If $\fg$ is the Lie algebra of infinitesimal coordinate changes (in our case $\fg = \text{Der}(\ah)$) then for any open contractible chart $U$ the action of $\fg (U)$ on $\call ^* (U)$ is homotopically trivial; moreover, for $X$ from $\fg$, let $L_X : \call ^* \rightarrow \call^*$ be the action of $\fg$ on $\call ^*;$ one can also construct operators $i_X :  \call ^* \rightarrow \call^{*-1}$ satisfying the usual formulas:
$$ [L_X, i_Y] = i_{[X,Y]} ; [i_X, i_Y] = 0 ; [d, i_X] = L_X $$
In our case $X = \frac{1}{i \hb} \text{ad}(a)$ and $i_X = L_{\underline{a}}$ where $\underline{a}$ denotes $a$ viewed as a Hochschild zero cochain.

Let $\widehat{\ah}(\R ^{2n})$ be the Weyl algebra $\C[[\varksi_1 , \ldots , \varksi _n, x_1 , \ldots , x_n , \hb ]]$ with the Moyal product. Let $\widehat{\fg} = \text{Der}(\widehat{\ah}(\R ^{2n}))$; then $\frak{u}(n)$ is a Lie subalgebra of $\widehat{\fg}.$ Put 
$$
    \ll ^* = \rc _{-*} (\widehat{\ah} (\R ^{2n}))[\hb ^{-1}]
$$

Let $W^*(\frak{u}(n))$ be the Weil algebra. As a partial case of a general construction one gets a cochain map
$$
   [W^*(\frak{u}(n)) \otimes \ll ^*]_{basic} \rightarrow \Cc (M, \call ^*)
\tag6.9
$$
The subscript $basic$ stands for $\{\alpha : L_X \alpha = i_X \alpha = 0, X \in \frak{u}(n) \}.$

To finish the proof of the theorem one constructs the fundamental class $\Omega$ in $ [W^*(\frak{u}(n)) \otimes \ll ^*]_{basic}$ and proves an analogue of Theorem 6.2 by an explicit calculation.  \pf
\enddemo

Let us finish by saying a few words about how to carry out all the above computations explicitly at the level of cochains. Recall that one can realize any deformation $\ah (M)$ as the space of horizontal sections of a Fedosov connection $\Nabla : \Omega ^0 (M, \W) \rightarrow  \Omega ^1 (M, \W)$ where $\W$ is the Weyl bundle (\c{Fe}). Then, instead of the Cech complex $\Cc ^* (M, \call ^*),$ one considers the de Rham complex $\Omega ^* (M, \ll ^*);$ all the cochains which participate in the proof have an easy explicit Chern - Weil style realization in this complex. 

\vskip .2in
\noindent\ub{\bf Section 7.}\quad\ub{\bf Examples}
\vskip .2in
For any differential graded algebra $B$ there is a natural gradind on
the Hochschild chain complex $C_*(B,B).$ We denote by  $C_*(B,B)^j$
the space of homogeneous elements of degree $j$ in  $C_*(B,B).$Let
$$
\gathered
\Cal{C}^n(A)=\bigoplus_{i+j=n} C_{-i}(\Cal{E}_A^*,\Cal{E}_A^*)^j \\
\Cal{C}^n_\infty(A)=\prod_{i+j=n} C_{-i}(\Cal{E}_A^*,\Cal{E}_A^*)^j
\endgathered
$$
We have defined a homotopically associative $\bullet$ product on these complexes; their homology is denoted by $\Cal{H}^*(A)$, resp. by $\Cal{H}^*_{\infty} (A)$.
\proclaim
{\ub{Theorem 7.1}} Let $A = k[X]$ be the algebra of regular functions on an affine algebraic variety $X$ over a field $k$ of characteristic zero. Then
$$
\Cal{H}^*(A) \rightisoarrow \Cal{D}(\Omega^*_{X/k})
$$
(the ring of all differential operators on $\Omega^*_{X/k}$);
$$
\Cal{H}^*_{\infty}(A) \rightisoarrow \text{End}(\Omega^*_{X/k})
$$
\endproclaim
Note that $HH_*(A) \simeq \Omega^*_{X/k}$; the $\bullet$ action of $\Cal{H}^*$ on $HH^*$ is the obvious one.
\demo
{\ub{Proof of 7.1}} 
\proclaim
{\ub{Lemma 7.2}} Let $\Cal{E}^*,$ $\Cal{F}^*$ be two differential graded algebras together with a quasi-isomorphism $f: \Cal{E}^* \to \Cal{F}^*.$ Then $f$ induces a quasi-isomorphism
$$
\bigoplus_{i+j=n} C_{-i}(\Cal{E}^*,\Cal{E}^*)^j \rightisoarrow \bigoplus_{i+j=n} C_{-i}(\Cal{F}^*,\Cal{F}^*)^j
$$
\endproclaim
\demo
{\ub{Proof}} Consider a filtration $F^p = \bigoplus_{i+j=n; i \geq p} C_{-i}(\Cal{E}^*,\Cal{E}^*)^j$ and a similar filtration for $\Cal{F}^*$. The map induced by $f$ is an isomorphism at the level of
$E^1$ terms of the spectral sequences with $E^1 = \bigoplus_{i+j=n}
C_{-i}(H^*(\Cal{E}^*),H^*(\Cal{E}^*))^j$, 
resp.

 $\bigoplus_{i+j=n}
C_{-i} (H^*(\Cal{F}^*),H^*(\Cal{F}^*))^j$, the Hochschild complex of
the graded algebra of cohomology of $\Cal{E}^*$, resp. $\Cal{F}^*.$ Note that we are able
to apply this argument to those converging spectral sequences (this is not true if one replaces $\bigoplus$ by $\prod$).     \pf
\enddemo
Note that by \c{HRK} $\E$ is quasi-isomorphic to $\Cal{F}^* = \Gamma (X, \bigwedge ^*TX)$ as differential graded algebras. If $x$ is a point of $X,$  for  germs at $x$ one has $\Cal{F}^*_x = \Cal{O}_x \otimes \bigwedge ^*(\partial_ {x_1}, \ldots, \partial_ {x_n});$ 
$HH_*(\Cal{F}^*)_x = \Omega^*_{X/k,x} \otimes \bigwedge ^*(\partial_{x_1}, \ldots, \partial_ {x_n}) \otimes S^*(d\partial_{x_1}, \ldots, d\partial_ {x_n})$ (classes $d\partial_{x_i}$ are represented by cycles $(1, d\partial_{x_i})$). Now, looking at the $\bullet$ action on $HH^*(\Cal{O}_x) \simeq \Omega ^*_{X/k,x}$, one sees that $\Omega ^*_{X/k, x}$ acts by multiplication, $\partial _{x_i}$ acts by contraction $i_{\partial _{x_i}}$ and $d\partial _{x_i}$ acts by Lie derivative $L_{\partial _{x_i}}.$ Whence the isomorphism
$$
\Cal{H}^*(\Cal{O}_x) \rightisoarrow \Cal{D}(\Omega^*_{X/k, x});
$$
this statement may be globalized using a standard argument.

Now let us compute $\Cal{H}^*_{\infty}(A)$. We have to use the other
spectral sequence whose $E_1$ term is $HH_*(\E)$ , the Hochschild
homology of the graded algebra $\E$ (if one forgets about the
differential $\delta$). 

For the  vector spaces $V$ and $W$ we will denote by $T(V')
\widehat{\otimes} W$ the space of multi-linear maps from $V$ to $W$. We
will use similar notation for certain subquotients of $T(V').$

  As a graded algebra, $\E$ is just $T^*(A') \widehat{\otimes} A$ by which we mean the space of all multilinear maps from $A$ to $A$. One can show that 
$$
   E_2^{-p, n} = HC^{n-1}(A)  \widehat{\otimes} \Omega ^p +  HC^{n-2}(A)  \widehat{\otimes} \Omega ^{p+1}
$$
$$
   E_3^{-p, n} = (\text{Ker}S) ^{n-1}  \widehat{\otimes} \Omega ^p +  (\text{Coker}S)^{n-2}  \widehat{\otimes} \Omega ^{p+1}
$$
where $S : HC^n \to HC^{n+2}$ is the Bott operator (\c{L}).
Recall the Gysin exact sequence (\c{C}, \c{T}, \c{L}):
$$
    0 \longrightarrow (\text{Coker}S)^{n} \longrightarrow HH^n(A) \longrightarrow (\text{ker}S)^{n-1} \longrightarrow 0 
$$
It is not hard to construct for any element of $HH^*(A) \widehat{\otimes} \Omega ^*_{X/k}$ a corresponding cycle of the complex $\Cal{C}^*$. This shows that the spectral sequence degenerates at $E_3$ term and that $\Cal{H}^*_{\infty}(A) \rightisoarrow HH^*(A) \widehat{\otimes} \Omega ^*_{X/k}$     \pf
\enddemo

\vskip .2in
\noindent\ub{\bf Section 8.}\quad\ub{\bf Characteristic cochain of a
flat element.}
\vskip .2in

Let $\n$ be an element of $C^*\asa$ \st $|\n|$ is odd (and therefore
$\deg\n$ is even).  Assume that
$$
        \delta\n+\fr12[\n,\n]=0
\tag8.1
$$
Consider an element 

$$
        e^{\n}=\sum\fr{1}{n!}\wedge^n\n\in\wedge^*\fg(A)
$$
in the completed space $\prod R_k({\frak{g}} (A))$.
Put
$$
        \X(\n)=J(e^{\n}\otimes 1)
\tag8.2
$$
This is a well-defined operator 
$$
        \X(\n) :CC_*^{\per,(0)}(A)\to CC_*^{\per}(A)
$$
where $CC_*^{\per,(0)}$ is the subcomplex of finite cochains:
$$
 CC_n^{\per, (0)}(A)=\underset{i\equiv n(2)}\to\bigoplus
        A\otimes\oa^{\otimes n},
$$
while, as usually,
$$
 CC_n^{\per}(A)=\underset{i\equiv n(2)}\to\prod
        A\otimes\oa^{\otimes n}.
$$
\proclaim
{\ub{Theorem 8.1}} 
$$
        [B+b,\;\X(\n)]=\X(\n)\cdot L_{\n}.
$$
\endproclaim

\demo
{\ub{Proof}}  Indeed, in the completion of $\wedge^*\fg(A)\otimes U(\fg(A))$:
$$
\gathered
        (\delta+\partial^{\text{Lie}})(e^{\n}\otimes 1)=
        e^{\n}(\delta\n+\fr12[\n,\n])\otimes   \\
        \otimes 1 + e^{\n}\otimes L_{\n}.  \qquad\text{\pf}
\endgathered
$$
\enddemo

In some cases there is a natural topology on $A$ such that $\X(\n)$
converges and defines an operator on a complex larger than $CC^{per,
0}$. 

The formula for $\X(\n)$ is obtained as follows.  In the 
periodic cyclic complex of the algebra generated by the even
elements $\n$ and $e^{-t\n}$, $t\geq 0$, put
$$
\gathered
        e^{I_{\unn}}=\underset{n\geq 0}\to\sum\;\;
        \underset{\Sb t_0+\dots+t_n=1 \\ t_i\geq 0 \endSb}\to\int
        e^{t_o i_{\n}}S_{\unn}e^{t_1 i_{\n}}\dots S_{\unn}
        e^{t_n i_{\n}}dt\dots dt_{n-1}  \\
        e^{ti_{\unn}}(D_0,\dots,D_m)=(D_0\smile e^{t\n},D_1,\dots,D_m);  \\
        S_{\unn}(D_0,\dots,D_m)=\overset{n}\to{\underset{j=0}\to\sum}\;
        \overset{j-1}\to{\underset{i=0}\to\sum}
        (-1)^{jn+n+j+i}(1,D_j,\dots,D_m,\;  \\
        D_0,\dots,D_i,\;\n,\dots,D_{j-1})
\endgathered
$$
Then
$$
        \X(\n)(\alpha) = \alpha \bullet (e^{I_{\unn}}1)
$$
where $\bullet$ is the product from Theorem 3.1.  One has 
$$
\gathered
        e^{I_{\unn}}1=e^{\n}-\int (e^{t_0\n},e^{t_1\n},\n)+  \\
        +\int[(2(e^{t_0\n},e^{t_1\n},\n,e^{t_2\n},\n)-  \\
        -(e^{t_0\n},\n,e^{t_1\n},\n,e^{t_2\n})+(e^{t_0\n},\n,
        e^{t_1\n},e^{t_2\n},\n)]+\dots
\endgathered
$$

\vskip .2in
\noindent\ub{\bf Section 9.}\quad\ub{\bf Bivariant JLO cochain.}
\vskip .2in
In \c{JLO} and \c{GS} a Chern character of a $\theta$-summable
Fredholm module over a Banach algebra $A$ was defined. It is given by
a cocycle of the complex dual to $CC^{per,0}(A)$; this cocycle
satisfies a special growth condition which makes it an entire cyclic
cocycle in the sense of Connes. Here we will construct its bivariant
version, which means, a homomorphism
from $CC^{per,0}(A)$ to $CC^{per}(A)$ whose composition with a trace is
given by the formula from  \c{JLO} and \c{GS}.

Given a $\Z_2$-graded algebra $A$ and an odd element $\sd$ of $A$,
put 
$$
        \n=ad(\sd)-\sd^2;
\tag9.1
$$
then $\n$ satisfies (8.1).

\proclaim
{\ub{Theorem 9.1}} Put
$$
        ch(\sd)=\X(\n)\cdot e^{L_{\ssd}};
$$
then
$$
        [B+b,ch(\sd)]=0
\tag9.2
$$
\endproclaim

\demo
{\ub{Proof}}  We have seen that
$$
        [B+b,\;\X(\n)]=\X(\n)L_{\n};
\tag9.3
$$
it suffices to show that 
$$
        [B+b,e^{L_{\underline{\ssd}}}]=-L_{\n}\cdot e^{L_{\underline{\ssd}}}.
\tag9.4
$$
Indeed,
$$
\gathered
        [B+b,e^{L_{\usd}}]=-\int_0^1 e^{tL_{\usd}}L_{ad\ssd}
        e^{(1-t)L_{\ssd}}dt=  \\
        =-L_{ad\ssd}\cdot e^{L_{\ssd}}+\int t\,L_{[\underline{\ssd\ssd}]}
        e^{L_{\ssd}}dt=  \\
        =-L_{\n}e^{L_{\usd}},  \qquad \text{\pf}
\endgathered
$$
\enddemo

\vskip .2in
\noindent\ub{\bf Section 10.}\quad\ub{\bf ``The Fukaya category''.}
\vskip .2in

In this Section we extend the results of Section 2 by ``moving away from the diagonal''. For any automorphism $\alpha$ of an algebra $A$ we define the twisted Hochschild complex $C_*(A,A_{\alpha})$ (related to noncommutative geometry of fixed points of $\alpha$) and construct in these terms a homotopically associative category for which $\Cal{C}^*(A)$ or $\Cal{C}^*_{\infty}(A)$ is the ring of endomorphisms of the identity.

Let $\alpha,$ $\beta$ be two automorphisms of $A.$ Define the new bimodule $_{\alpha} A _{\beta}$ over $A$ as follows: $ _{\alpha} A _{\beta} = A$ as $k$-modules and $a\cdot m \cdot b = \alpha(a) m \beta (b)$ for $a,m,b$ in $A.$
Consider the chain complex $C_*(A,\;_{\alpha} A _{\beta})$ as in Section 2 and the cochain complex $_{\alpha}\Cal{E}^*_{\beta} = C^*(A,\;_{\alpha} A _{\beta})$ as in Section 1. We shall write $A_{\alpha} = _{\text{id}} A _{\alpha}$. The cup product on \hoc cochains is well defined as a morphism $_{\alpha}\Cal{E}^*_{\beta} \otimes _{\beta}\Cal{E}^*_{\gamma} \rightarrow_{\alpha}\Cal{E}^*_{\gamma}$. 

Given an automorphism $\alpha$ of $A$, one can define its action on a \hoc 
cochain from $\E$ in two ways:

$$
(\alpha D)(a_1, \ldots, a_n) = \alpha (D(a_1, \ldots, a_n));
$$
$$ 
(D \alpha)(a_1, \ldots, a_n) =  D(\alpha a_1, \ldots,\alpha a_n)
$$
Both yield morphisms of complexes $\E \rightarrow _{\alpha}\Cal{E}_{\alpha}$.
We define an $\E$-bimodule structure on $_{\alpha}\Cal{E}^*_{\beta} $ as follows:
$$
D\cdot M \cdot E = D\alpha \smile M \smile \beta E 
$$
for $D,E \in \E$ and $M \in $$_{\alpha}\Cal{E}^*_{\beta} $.

\proclaim
{\ub {Theorem 10.1}}
The pairings (0.6) can be extended to the homotopically associative natural morphisms of complexes
$$ \bullet :C_*(A,A_{\alpha}) \otimes C_*(\E, \;_{\alpha}\Cal{E}^*_{\beta}) \rightarrow  C_*(A, A_{\beta})
\tag10.1
$$
$$
\bullet :C_*(\E,\; _{\alpha}\Cal{E}^*_{\beta}) \otimes C_*(\E, \;_{\beta}\Cal{E}^*_{\gamma}) \rightarrow C_*(\E, \;_{\alpha}\Cal{E}^*_{\gamma})
\tag10.2
$$
\endproclaim
{\ub {Proof}}
As in Section 2, put for $a \in C_*(A,A_{\alpha})$ and $x \in C_*(\E, \; _{\alpha}\Cal{E}_{\beta})$
$$         
     a\bullet x =  a\bullet _1 x +  a\bullet _2 x
$$
where
$$
  \multline
     (a_0 ,\ldots , a_n) \bullet_1  (D_0 , \ldots ,      D_m) = \\
\sum \pm 
 (a_0 \cdot D_0(a_1, \ldots ) ,\ldots ,a_{i_1} , D_1 (a_{i_1 + 1}, \ldots ),\ldots ,D_m (a_{i_m + 1}, \ldots ) , \ldots)
\endmultline
\tag10.3
$$
$$
 \multline
     (a_0 ,\ldots ,a_n) \bullet_2  (D_0 , \ldots ,      D_m) = \\
=\sum_{q \leq n+1} \pm (D_m (a_q, \ldots , a_n, a_0, \alpha a_1, \ldots, \alpha a_{i_0})\cdot D_0(a_{i_0 + 1}, \ldots ) ,\ldots , A_{i_1} , \\
, D_1 (a_{i_1 + 1}, \ldots ), \ldots , D_{m-1} (a_{i_{m-1} + 1}, \ldots ) , \ldots )
\endmultline
\tag10.4
$$

For $A \in  C_*(\E, \; _{\alpha}\Cal{E}_{\beta})$ and $x \in  C_*(\E, \; _{\beta}\Cal{E}_{\gamma})$

$$         
     A\bullet x =  A\bullet _1 x +  A\bullet _2 x
\tag10.5
$$
where
$$
  \multline
     (A_0 ,\ldots , A_n) \bullet_1  (D_0 , \ldots ,      D_m) = \\
\sum \pm 
 (A_0 \smile D_0\{A_1, A_2, \ldots \} ,\ldots ,A_{i_1} , D_1 \{A_{i_1 + 1}, \ldots \},\ldots ,D_m \{A_{i_m + 1}, \ldots \} , \ldots)
\endmultline
\tag10.6
$$
$$
 \multline
     (A_0 ,\ldots ,A_n) \bullet_2  (D_0 , \ldots ,      D_m) = \\
=\sum_{q \leq n+1} \pm (D_m\{\{A_q, \ldots , A_n, A_0, A_1, \ldots, A_{i_0}\}\} \smile D_0\{A_{i_0 + 1}, \ldots \} ,\ldots , A_{i_1} , \\
, D_1 \{A_{i_1 + 1}, \ldots \}, \ldots , D_{m-1} \{A_{i_{m-1} + 1}, \ldots \} , \ldots )
\endmultline
\tag10.7
$$
and 
$$
\multline
D_m\{\{A_q, \ldots , A_n, A_0, A_1, \ldots, A_{i_0}\}\}(a_1, \ldots, a_n) = \\
= \sum \pm D_m (\alpha a_1, \dots, A_q (\alpha a_j, \dots, ), \dots,\\
 \alpha a_t,  
 A_0 (a_{t+1}, \dots, ),\beta a_p, \dots, \beta A_1(a_r, \dots ), \beta a_s, \dots )
\endmultline
$$

One checks that these maps are homotopically associative morphisms of complexes.
\pf

One gets the cohomology groups $\Cal{H}^*(\alpha, \beta)$ and  $\Cal{H}^*_{\infty}(\alpha, \beta)$ which form a category. When $A$ is the ring of functions on $X$ then, in simple cases, the groups $\Cal{H}^*(\alpha, \beta)$ are results of some standard $D$-module constructions on $T^*(X).$ When $A$ is a deformed algebra of functions on a symplectic manifold $M$ then , when $M$ is simply connected, all the automorphisms of $A$ are of the form $exp(\text{ad} \frac{1}{i\hb}f);$ $\Cal{H}^*(\alpha, \beta)$ and  $\Cal{H}^*_{\infty}(\alpha, \beta)$ reflect the geometry of fixed point sets of $\alpha$ and $\beta$ and of related loop and path spaces. It looks as if there was, along with $\Cal{H}^*$ and  $\Cal{H}^*_{\infty},$ the intermediate semi-infinite cohomology which is more closely related to Floer cohomology of loop spaces.


\Refs\nofrills{BIBLIOGRAPHY}

\widestnumber\key{NT1}
\ref
  \key B
  \by J.L. Brylinski
  \paper Some examples of Hochschild and cyclic homology
  \jour Lecture Notes in Math.
  \vol 1271
  \yr 1987
  \pages 33-72
  \endref

\ref
  \key BFFLS
  \by  F. Bayen, M. Flato, C. Fronsdal, A. Lichnerowiscz, D. Sternheimer
%  \paper  
1  \jour  Ann. Phys.
  \vol  \underbar{111}
  \yr  1978
  \pages  61--151
  \endref

\ref
  \key BG
  \by J.L Brylinski and E.Getzler 
  \paper The homology of algebras of pseudo-differenial symbols and the noncommutative residue
  \jour K-Theory
  \vol 1
  \yr 1987
  \pages 385-403
  \endref


\ref
  \key CE
  \by A. Cartan and S. Eilenberg
  \paper Homological algebra
  \paperinfo Princeton University Press, 1956
\endref
\ref
  \key C
  \by A.Connes
  \paper Noncommutative differential geometry
  \paperinfo Publ. Math. IHES, 1986
\endref

\ref
  \key  CFS
  \by  A. Connes, M. Flato, D. Sternheimer
  \paperinfo  Closed star products and cyclic 
homology, preprint, Paris, 1991
 % \vol  
 % \yr  
 % \pages  
  \endref
\ref
  \key DT
  \by Y. Daletsky and B. Tsygan
  \paper Operations on Hochschild and cyclic complexes
  \paperinfo preprint 1992
\endref
 \ref
  \key De
  \by M. Deligne
  \paper  D\'{e}formations d'alg\`{e}bre des fonctions
d'une variet\'{e} symplectique: comparaison entre Fedosov et DeWilde,Lecomte
  \jour Selecta Math.
  \vol 1
  \yr 1995
  \pages 667-697        
\endref         
 \ref
  \key DWL
  \by M. De Wilde and P.B.A.Lecomte
  \paper  Existence of star-products and of formal deformations in Poisson Lie algebras of arbitrary Poisson manifolds
  \jour Lett. Math. Phys.
  \vol 7
  \yr 1983
  \pages 487-496        
\endref
\ref
  \key Fe
  \by B. Fedosov
  \paper Index theorems
  \paperinfo Itogi Nauki i Tekhniki, vol. 65, 1991 (in Russian) 
\endref
\ref
  \key FT
  \by B. Feigin and B.Tsygan
  \paper Additive K-theory
  \jour Lecture Notes in Math.
  \vol 1289
  \yr 1987
  \pages 66-220         
\endref
\ref
  \key Fu
  \by K. Fukaya
  \paper Morse homotopy, $A_{\infty}$ categories and Floer homologies
  \paperinfo MSRI preprint No. 020-94, 1993.
\endref
\ref
  \key Ge
  \by M. Gerstenhaber
  \paper The cohomology structure of an associative ring
  \jour Ann. of Math.
  \vol 78 
  \yr 1963
  \pages 267--288       
\endref                 
\ref
  \key GV
  \by M. Gerstenhaber and A. Voronov
  \paper Homotopy $G$-algebras and moduli space operad
  \jour Int. Math. Res. Notices
  \yr 1995
  \pages 141-153        
\endref

\ref
  \key G
  \by E. Getzler
  \paper Cartan homotopy formulas and the Gauss-Manin connection
         in cyclic homology
  \jour Israel Math. Conf. Proc.
  \vol 7
  \yr 1993
  \pages 65--78         
\endref

\ref
  \key GJ
  \by E.Getzler and J.Jones
  \paper Operads, homotopy algebra and iterated integrals for double 
loop spaces
  \paperinfo preprint hep-th/9403055, 1994
  \endref

\ref
  \key GJ1
  \by E.Getzler and J.Jones
  \paper $A_{\infty}$ algebras and the cyclic bar complex
  \jour Illinois J. of Math.
  \vol 34
  \yr 1990
  \pages 256-283
  \endref

\ref
  \key GJP
  \by E.Getzler, J.Jones and S.Petrack
  \paper Differential forms on loop spaces and the cyclic bar complex
  \jour Topology
  \vol 30
  \yr 1991
  \pages 339-371
  \endref

 \ref
  \key GS
  \by E.Getzler and A.Szenes
  \paper On the Chern character of a theta-summable Fredholm module
  \jour J. Funct. Anal.
  \vol 84
  \yr 1989
  \pages 343-357        
\endref 

\ref
  \key HJ
  \by C.E.Hood and J.Jones
  \paper Some algebraic properties of cyclic homology groups  
\jour K-theory
\vol 1
\pages 361-384
  \endref
\ref
  \key HRK
  \by G.Hochschild, B.Kostant and A.Rosenberg
  \paper Differential forms on regular affine algebras
  \jour Transactions AMS
  \vol 102
  \yr 1962
  \pages 383-408                
\endref
\ref
  \key JLO
  \by A. Jaffe, A. Lesniewski and K. Osterwalder
  \paper Quantum K-theory I
  \jour K-theory
  \vol 2
  \yr 1989
  \pages 675-682                
\endref
\ref
  \key LQ
  \by J.L.Loday and D. Quillen
  \paper Cyclic homology and Lie algebra homology of matrices
  \jour Comment. Math. Helv.
  \vol 59
  \yr 1984
  \pages 565-591                
\endref
\ref
  \key M  \by Yu. I. Manin
  \paper Topics in non-commutative geometry
  \paperinfo Princeton University Press, 1991
\endref
\ref
  \key L
  \by J.L.Loday
  \paper Cyclic Homology
  \jour Springer Verlag
  \vol 
  \yr 1993
\endref
\ref
   \key MS
   \by D. MacDuff and D. Salomon
   \paper Quantum Comology and J-holomorphic Curves
   \paperinfo AMS University Lecture Series
  \yr 1994
\endref
\ref
  \key NT
  \by R. Nest and B. Tsygan
  \paper Algebraic index theorem for families
  \jour Advances in Math
  \vol 113
  \pages 151-205
  \yr 1995
\endref
 \ref
  \key NT2
  \by R. Nest and B. Tsygan
  \paper Homological properties of the category of algebras  and       characteristic classes
  \paperinfo Preprint, Univ. of Heidelberg, 1995
\endref
\ref
  \key R
  \by G. Rinehart
  \paper Differential forms on commutative algebras
  \jour Trans. AMS
  \vol 108
 \yr 1963
  \pages 139-174
\endref
\ref
  \key T
  \by  B. Tsygan
  \paper  Homology of matrix Lie algebras over 
rings and Hochschild homology
  \jour  Uspekhi Mat. Nauk
  \vol  38, 2
  \yr  1983
  \pages  217-218
 \endref
\ref
  \key W1
  \by  M. Wodzicki
  \paper  Cyclic homology of differential operators
  \jour  Duke Mathematical Journal
  \vol  54
  \yr  1987
  \pages  641-647
 \endref
\ref
  \key W2
  \by  M. Wodzicki
  \paper  Cyclic homology of pseudodifferential operators and noncommutative Euler class
  \jour  C.R.A.S.
  \vol  306
  \yr  1988
  \pages  321-325
 \endref
\endRefs
\enddocument